\documentclass[reqno]{amsart}
 \usepackage{amsbsy,amssymb,amscd,amsfonts,latexsym,amstext,delarray,
 amsmath,graphicx, color}
 \usepackage[dvips]{epsfig}
\usepackage{hyperref}
\input xypic

\definecolor{Green}{rgb}{0.0,0.40,0.0}

\newtheorem{thm}{Theorem}[section]
\newtheorem{prop}[thm]{Proposition}
\newtheorem{cor}[thm]{Corollary}
\newtheorem{lem}[thm]{Lemma}

\newtheorem{defn}[thm]{Definition}
\newtheorem{rem}[thm]{Remark}

\def\qqq{\,,\quad~\forall}

\def\Aut{{\rm Aut}}

\def\Hom{{\rm Hom}}

\def\Sp{{\rm Spec}}

\def\A{{\mathbb A}}
\def\C{{\mathbb C}}
\def\F{{\mathbb F}}
\def\K{{\mathbb K}}

\def\N{{\mathbb N}}

\def\Q{{\mathbb Q}}
\def\R{{\mathbb R}}
\def\Z{{\mathbb Z}}

\def\cA{{\mathcal A}}
\def\cB{{\mathcal B}}
\def\cC{{\mathcal C}}

\def\cE{{\mathcal E}}

\def\cG{{\mathcal G}}
\def\cH{{\mathcal H}}

\def\cF{{\mathcal F}}

\def\cO{{\mathcal O}}
\def\cP{{\mathcal P}}

\def\cR{{\mathcal R}}

\def\cT{{\mathcal T}}

\newcommand{\ie}{{\it i.e.\/}\ }
\newcommand{\eg}{{\it e.g.\/}\ }
\newcommand{\cf}{{\it cf.\/}\ }

\def\ker{{\mbox{Ker}}}

\def\Hom {{\mbox{Hom}}}

\def\Gm{{\mathbb G}_m}
\def\truc{\cF}

\newcommand{\nil}[1]{}

\title
{Fun with $\F_1$}

\author[Connes]{Alain Connes}
\author[Consani]{Caterina Consani}
\author[Marcolli]{Matilde Marcolli}
\address{A.~Connes: Coll\`ege de France \\
3, rue d'Ulm \\ Paris, F-75005 France
\\ I.H.E.S. and Vanderbilt
University} \email{alain\@@connes.org}
\address{C.~Consani: Mathematics Department \\ Johns Hopkins
University \\ Baltimore, MD 21218 USA} \email{kc\@@math.jhu.edu}
\address{M.~Marcolli: Max--Planck Institut f\"ur Mathematik  \\
Vivatsgasse 7 \\
Bonn, D-53111 Germany} \email{marcolli\@@mpim-bonn.mpg.de}

\date{}
\begin{document}
\maketitle \vspace{2cm}

\begin{abstract} We show that the algebra and the
endomotive of the quantum statistical mechanical system of Bost--Connes
naturally arises by extension of scalars from the ``field with one element" to
rational numbers. The inductive structure of the abelian part of the endomotive
corresponds to the tower of finite extensions of that ``field", while the
endomorphisms  reflect the Frobenius correspondences. This gives in particular
an explicit model over the integers for this endomotive, which is related to
the original Hecke algebra description. We study the reduction at a prime of
the endomotive and of the corresponding noncommutative crossed product algebra.
\end{abstract}

 \tableofcontents

\section{Introduction}

Starting with seminal observations of J.~Tits on the classification of simple
finite groups (\cf~\cite{Tits}), the a priori vague idea that a suitable
analogue of the geometry over the finite fields $\F_q$ should make sense in the
limit case ``$q=1$'' has been taking more and more substance and has given rise
to a number of different approaches (\cf \cite{Kapranov}, \cite{Manin},
\cite{Soule99}, \cite{Soule}, \cite{Durov}, \cite{Haran}, \cite{TV}). So far,
the relation between these constructions and the Riemann zeta function has
remained elusive, in spite of the hope of being able to consider the tensor
product $\Z\otimes_{\F_1}\Z$ as a non-trivial analogue of the product of a
curve by itself (see \cite{Manin}).

It is known that the quantum statistical mechanical system of  \cite{BC} (which
we refer to as the BC-system) gives, after passing to the dual system, a
spectral realization of the zeros of the Riemann zeta function, as well as a
trace formula interpretation of the Riemann-Weil explicit formulas (see
\cite{BC}, \cite{Co-zeta}, \cite{CCM}, \cite{CCM2}, \cite{Meyer}).

The main result of the present paper is that the BC-system and the associated
algebraic endomotive as defined in \cite{CCM} appear from first principles, by
studying the algebraic extensions of $\F_1$ and implementing the techniques
developed in \cite{Kapranov} and \cite{Soule}.

In this formalism, a variety (of finite type) over $\F_1$ determines, after
extension of scalars to $\Z$, a variety over $\Z$. Moreover, even though the
algebraic nature of $\F_1$ is still mysterious, a basic equation of the theory
is the formal equality
\begin{equation}\label{basic}
\F_{1^n}\otimes_{\F_1} \Z:=\Z[T]/(T^n-1)\,, \qquad n\in\N .
\end{equation}
Our starting point is that the natural inductive structure defined by the
extensions $\F_{1^n}\subset \F_{1^m}$, for $n|m$, translates into a natural
inductive system of algebras, whose limit is the group ring $\Z[\Q/\Z]$. After
tensoring by $\Q$, this group ring is a key ingredient in the definition of the
BC-endomotive \cite{CCM}, since it describes the abelian part of the structure.
The second key ingredient is the semigroup of endomorphisms of the above
algebra associated to the action of $\N$, given by multiplication, in the group
$\Q/\Z$. This operation describes an analogue over $\F_1$ of the classical
Frobenius correspondence, and part of our investigation is directed at making
this statement more precise.

In \S \ref{BCrecall} we start our study by recalling the group theoretical and an
 equivalent geometrical description of the abelian part of the BC algebra.
 In particular, the sub-section \S \ref{BCrecall1} introduces the relevant abelian
algebra, that is the group ring of the abelian group $\Q/\Z$, together with the
endomorphisms $\sigma_n$ given by multiplication by $n$ in $\Q/\Z$ as well as
their partial inverses $\rho_n$. In \S \ref{alggeom} we describe the same space
by using elementary techniques of algebraic geometry.

In \S \ref{F1roots1} we give an interpretation of the abelian part of the
BC-system  in terms of a system of pro-affine varieties $\mu^{(k)}$ over $\F_1$ which are defined by considering affine group
schemes of roots of unity.

 The core of the paper concentrates on the definition
 of an integral model for the BC-system over $\Q$ and
 on the generalization of the notion of endomotive that
  was originally developed over fields of characteristic
  zero only. This study is motivated by the idea
  to achieve an interesting link between the
  thermodynamical system associated to the
  BC-algebra (and its connection to the
  zeta-function) and the theory of algebraic varieties over $\F_1$. Working with spaces over $\F_1$ implicitly requires one to define a geometric theory over $\Z$ and thus to set-up a corresponding construction over finite fields (and their extensions), after taking the reduction at the various primes. The main result of \S \ref{BCintegral} is that the original noncommutative BC-algebra $\cA_\Q$ has a model over $\Z$.

 In \S \ref{BCcstar} we shortly review the description of the $C^*$-algebra  of the
BC-system in terms of groupoids and in \S \ref{bcqsect} we recall the
presentation of $\cA_\Q$ by generators and relations. \S \ref{secttilderho}
describes how to eliminate the denominators in the partial inverses $\rho_n$ of
the endomorphisms $\sigma_n$. This leads us in \S \ref{bczz} to the definition
of the integral model $\cA_\Z$ of the BC-algebra by generators and relations.
In Proposition \ref{pres4} and Corollary \ref{corpres4}, we show that the
general element of $\cA_\Z$ can be uniquely written as a sum of simple
monomials labeled by $\Q/\Z\times \Q_+^*$. In \S \ref{bchecke} we define an
isomorphism of $\cA_\Z$ with the integral version $\cH_\Z(\Gamma,\Gamma_0)$ of
the original Hecke algebra of \cite{BC}, and deduce from that the existence of
two {\em different} involutions on the rational algebra $\cA_\Q$. In \S
\ref{sectcharp} we analyze the BC endomotive over a perfect field of
characteristic $p$, and we relate the endomorphisms of this system to the
Frobenius correspondence.

We isolate the $p$-part $\cC_p$ of the BC-algebra in characteristic $p$ and
exhibit its nilpotent nature by constructing in Proposition
\ref{proptriangular} a faithful representation of $\cC_p$ as lower triangular
infinite matrices. This representation is obtained by relating the algebra
$\cC_p$ to a sub-semigroup of the  group of affine transformations of the
additive group $S=\cup \frac{1}{p^n}\Z\subset \R$. A new feature that arises in
positive characteristic is the appearance of unreduced algebras in the abelian
part of the system. We explain the effect of reduction of these algebras in \S
\ref{sectreduc} and briefly discuss in \S \ref{sectendounred}, the required
extension of the notion of endomotive to the general (unreduced) framework.

Finally, in \S \ref{BCF1sec} we prove that the BC system has a model defined
over $\F_1$. This result allows us to deduce that the symmetries of the BC
system are recovered from the automorphisms of $\F_{1^\infty}=\varinjlim
\F_{1^n}$ over $\F_1$. In fact, we show that the BC endomotive embodies the
structure of the extensions $\F_{1^n}$ of $\F_1$ through the Frobenius
correspondence which is implemented by the action of the endomorphisms on the
abelian part of the associated algebra. More precisely, we show that these
endomorphisms coincide with the Frobenius correspondence in the reduction of
the BC system over a perfect field of positive characteristic. We then use this
result to prove that the original analytic endomotive of the BC system can be
recovered from the data supplied at infinity in the form of an inductive system
of Banach algebras.
\medskip

{\bf Acknowledgment.} The authors are partially supported by NSF grants
DMS-FRG-0652164, DMS-0652431, and DMS-0651925.
The second author gratefully thanks l'Institut
des Hautes \'Etudes Scientifiques for the hospitality,
the pleasant atmosphere and the support received during
 a visit in January-April 2008. The third author would like to
thank Abhijnan Rej for some useful conversations.

\medskip

\section{The abelian part of the BC system and its endomorphisms}\label{BCrecall}

In this section we shall give a short overview of two equivalent
formulations of the abelian part of the  algebra describing the quantum
statistical mechanical system introduced in \cite{BC} as well as the associated endomotive \cite{CCM}. In the following and
throughout the paper we shall refer to it as the BC-system (\cf~Definition~\ref{bcsystem}).

\subsection{Group theoretic description}\label{BCrecall1}

The BC-endomotive over $\Q$ is defined as the algebraic crossed product of the
group ring $\Q[\Q/\Z]$ by the action of a semigroup of endomorphisms.

In the following, we denote by $e(r)$, for $r\in \Q/\Z$, the canonical
generators of the group ring $\Q[\Q/\Z]$ with presentation
\begin{equation}\label{pres}
    e(a+b)=e(a)e(b)\qqq a, b \in \Q/\Z.
\end{equation}

We now describe in group theoretic terms the semigroup action on this group
ring. Let $\Gamma=\Q/\Z$. For each $n$, let $\Gamma_n\subset \Gamma$ be the
$n$-torsion subgroup
\begin{equation}\label{ntorsion}
    \Gamma_n=\{x\in \Gamma\,|\, nx=0\}\,.
\end{equation}

\begin{prop}\label{propendo} Let $n\in \N$.
\begin{enumerate}
  \item[(a)] One has an exact sequence of abelian groups
  \begin{equation}\label{absequ}
    1\to \Gamma_n\to \Gamma\stackrel{\times n}{\to} \Gamma\to 1 .
  \end{equation}
  \item[(b)] The operator
  \begin{equation}\label{idem}
    \pi_n=\frac 1n\sum_{s\in\Gamma_n} e(s)
  \end{equation}
  defines an idempotent $\pi_n\in \Q[\Q/\Z]$. One has $\pi_n\pi_m=\pi_k$ where
  $k$ is the least common multiple of $n$ and $m$.
\item[(c)] The formula
\begin{equation}\label{endo}
\rho_n: \Q[\Q/\Z]\to\Q[\Q/\Z],~ \rho_n(e(r))=\frac 1n\sum_{ns=r} e(s)
\end{equation}
defines an endomorphism of $\Q[\Q/\Z]$. Moreover $\rho_n$ is a  ring
isomorphism between $\Q[\Q/\Z]$ and the reduced algebra by $\pi_n$, more
precisely one has
\begin{equation}\label{endobis}
\rho_n: \Q[\Q/\Z]\stackrel{\sim}{\to}\pi_n\Q[\Q/\Z]\,, \   \rho_n(e(r))=\pi_n
e(s) \qqq s\,~\text{s.t.}~ \ ns=r\,.
\end{equation}
\item[(d)] The formula
\begin{equation}\label{endosigma}
\sigma_n: \Q[\Q/\Z] {\to}\Q[\Q/\Z],~ \sigma_n(e(r))= e(nr) \qqq r\in \Q/\Z\,,\
n\in\N
\end{equation}
defines an endomorphism of $\Q[\Q/\Z]$ and one has
\begin{equation}\label{rhosigma1}
\sigma_n \rho_n (x) = x ,
\end{equation}
\begin{equation}\label{rhosigma2}
\rho_n \sigma_n (x) = \pi_n x.
\end{equation}
\end{enumerate}
\end{prop}

\proof (a) follows from the divisibility of the group $\Gamma$, which implies
the surjectivity of the multiplication by $n$.

(b) One checks that $\pi_n^2=\pi_n$ using \eqref{pres}. Given integers $n$ and
$m$, with $\ell=(n,m)$ their gcd, the map $(s,t)\in \Gamma_n\times
\Gamma_m\mapsto s+t\in \Gamma_k$ is an $\ell$ to $1$ map onto $\Gamma_k$, where
$k=nm/\ell$ is the  least common multiple of $n$ and $m$. Thus
$\pi_n\pi_m=\pi_k$.

(c) First, the homomorphism $\rho_n$ is well defined since
$$
\pi_n e(s)=\frac 1n\sum_{nu=s} e(u)
$$
is independent of the choice of $s$ such that $ns=r$. It defines an algebra
homomorphism, since with $ns_j=r_j$ one has
$$
\rho_n(e(r_1+r_2))=\pi_n e(s_1+s_2)=\pi_n e(s_1)\pi_n
e(s_2)=\rho_n(e(r_1))\rho_n(e(r_2)).
$$

We then show that $\rho_n$ is an isomorphism with the reduced algebra. We let
$\hat\Gamma$  be the (Pontrjagin) dual of $\Gamma$, \ie the group  of
characters of the group $\Gamma$. We let $E_n$ be the open and closed subset of
$\hat\Gamma$ given by the condition $\chi\in E_n$ if and only if
$\chi(\pi_n)=1$. This holds if and only if $\chi|_{\Gamma_n}=1$. This allows
one to identify the closed subset $E_n\subset \hat\Gamma$  with the space of
characters of the quotient group $\Gamma/\Gamma_n$. Using the identification of
$\Gamma/\Gamma_n$ with $\Gamma$ determined by the isomorphism \eqref{absequ},
one gets an isomorphism of $E_n \subset \hat\Gamma$ with  $\hat\Gamma$, given
by
\begin{equation}\label{formofrhon}
E_n\stackrel{\rho_n^*}{\to}\hat\Gamma,\qquad\chi\mapsto
\rho_n^*(\chi),\qquad\rho_n^*(\chi)(r)=\chi(s)\qqq s\,~\text{s.t.}~ \ ns=r.
\end{equation}
In other words, $\rho_n^*(\chi)=\chi\circ\rho_n$ at the level of the group
ring. The range of the algebra homomorphism $\rho_n$  is contained in the
reduced algebra by $\pi_n$, and $\rho_n(f)=0$ implies $\langle
\rho_n^*(\chi),f\rangle=0$ for all $\chi \in E_n$. It then follows that $f=0$
since $\rho_n^*$ is surjective, hence $\rho_n$ is injective.

To show that the algebra homomorphism $\rho_n$ is surjective on the reduced
algebra by $\pi_n$, it is enough to show that the range contains the
$\pi_ne(a)$ for all $a\in \Q/\Z$. This follows from \eqref{endo}.

(d) The map $r\in \Q/\Z\mapsto nr \in \Q/\Z$ is a group homomorphism. One has
$\sigma_n(\pi_n)=1$ by \eqref{idem} and one gets \eqref{rhosigma1} using
\eqref{endo}. One checks \eqref{rhosigma2} on the generators $e(s)$ using
\eqref{endo}.
\endproof

\subsection{The endomorphisms $\rho_n$ from algebraic
geometry}\label{alggeom}

Let us first recall the geometric construction introduced in \cite{CCM} which
gives rise to interesting examples of algebraic endomotives. One lets $(Y,y_0)$
be a pointed smooth algebraic variety (over a field $\K$ of characteristic
zero) and $S$ an abelian semi-group of algebraic self-maps $s:Y\to Y$ with
$s(y_0)=y_0$, which are finite (of finite degree) and unramified over $y_0$. In
this way one then obtains:\medskip

$\bullet$~A projective system of algebraic varieties
\[
X_s=\{y\in Y\,|\,s(y)=y_0\}\,, \ \ \xi_{s,s'}: X_{s'} \to X_s \ \ \ \
\xi_{s,s'}(y)=r(y)~\text{if}~s'=rs.
\]\smallskip

$\bullet$~Algebraic morphisms
\begin{equation} \label{betas}
\beta_s: X=\varprojlim X_u \to X^{e_s}\,, \ \ \xi_u\beta_s(x)=s\xi_u(x)
\end{equation}
where $X^{e_s}=\xi_s^{-1}(y_0)\subset X$ are open and closed subsets of $X$.\medskip

In other words, one obtains in this way a first action $\sigma_s$ of the
semigroup $S$ on the projective limit $X=\varprojlim X_u$, since the maps
$\beta_s$ given by applying $s$ componentwise commute with the connecting maps
of the projective system,
\begin{equation}\label{sigmabeta}
\sigma_s(f)=f\circ \beta_s.
\end{equation}

In fact, since the maps $\beta_s$ are isomorphisms of $X$ with $X^{e_s}$, it is
possible to invert them and define a second action of $S$ that corresponds, at
the algebraic level, to the endomorphisms
\begin{equation}\label{defnendo}
\rho_s(f)(x) = \begin{cases} f(\beta_s^{-1}(x))&\text{if}~x\in X^{e_s}\\
0 &\text{if}~x\notin X^{e_s}.\end{cases}
\end{equation}

The BC-endomotive is a special case of this general construction. It
corresponds to the action of the semigroup $S=\N$ by finite algebraic
endomorphisms, unramified over $1$, on the pointed algebraic variety
$({\Gm}_{/\Q},1)$ (\cf~\cite{CCM} Proposition 3.7). One has
${\Gm}_{/\Q}=\rm{Spec}(\Q[T^{\pm 1}])$ and the action of $\N$ is given by
$a\mapsto a^n$ on the coordinate $a\in \bar \Q$ of the point associated to the
character $P(T)\mapsto P(a)$. Equivalently, at the algebra level, this action
is described by
\begin{equation}\label{geoendo}
\Q[T^{\pm 1}] \ni P(T,T^{-1})\mapsto P(T^n,T^{-n}) \in \Q[T^{\pm 1}].
\end{equation}
The point $1\in {\Gm}_{/\Q}$ is a fixed point and its inverse image under the algebraic
map $a\mapsto a^n$ is  $X_n = \rm{Spec}(\Q[T^{\pm 1}]/(T^n-1))$. The spaces
$X_n$ form a projective system indexed by $\N$, with partial order given by
divisibility. That is, for a pair of natural numbers $r,n$ with $r = ns$, we
have maps
\[
\xi_{n,r}: X_{r} \to X_n,\quad x\mapsto x^s \,.
\]
One lets $X=\varprojlim_n X_n$. The base point $1$ belongs to $X_n$ for all $n$
and defines a component $Z_n=\{1\}$ of $X_n$.  One checks that the description
of the algebra morphisms $\rho_n$ given by \eqref{defnendo} agrees with that
given in \eqref{endo}.  Here  the closed and open subset $X^{e_n}\subset X$ of
$X$ is simply the inverse image $\xi_{n}^{-1}(Z_n)\subset X$ of $Z_n$ by the
canonical map $\xi_{n}$ from the projective limit $X=\varprojlim_n X_n$ to
$X_n$.

The relation between this geometric description of the BC-endomotive and the
previous group theoretic one can be seen in the following way.

Let $u(n)$ be the class of $T$ mod. $T^n-1$, \ie  the canonical generator of
the algebra $\Q[T^{\pm 1}]/(T^n-1))$. Then the homomorphism $\xi_{m,n}$ is
given by
\begin{equation}\label{morphism}
    \xi_{m,n}(u(n))=u(m)^a\,, \ \ a=m/n\,.
\end{equation}
The isomorphism with the group theoretic description is then obtained by
mapping $u(n) \mapsto e(\frac 1n)\in \Q[\Q/\Z]$.

\section{$\F_{1^\infty}$ and the abelian part of the BC-system}\label{F1roots1}

In this section we describe the group ring part of the algebra of the BC-system
in terms of schemes of finite type over $\F_1$ in the sense of \cite{Soule}.
This is done by introducing a family of affine algebraic varieties $\mu^{(n)}$
over $\F_1$. We will show that these spaces can be organized in two ways: as an
inductive system related to the affine multiplicative group scheme over $\F_1$
and also as a pro-variety $\mu^{(\infty)} =\varprojlim_n \mu^{(n)}$. The
relation with the BC-system arises exactly when one works with the second
system.

We first recall the examples of extensions $\F_{1^n}$ of $\F_1$,  developed in
\cite{Kapranov} and \cite{Soule}, which are the analogues for $q=1$ of the
field extensions $\F_{q^n}$ of $\F_q$. The main idea is that these extensions
are described by the formula \eqref{basic} after extending the coefficients
from $\F_1$ to $\Z$. Notice that neither $\F_1$ nor its extensions $\F_{1^n}$
need to be properly defined for \eqref{basic} to make sense. Following
\cite{Kapranov}, while ``vector spaces over $\F_1$" correspond to sets, those
defined over the extension $\F_{1^n}$ correspond to sets with a free action of
the group $\Z/n\Z$.
\smallskip

When $n|m$, one expects in analogy with the case of $\F_{q^n}$ ($q>1$ a
rational prime power) to have extensions
\begin{equation}\label{fieldext}
    \F_{1^n}\subset \F_{1^m}\,
\end{equation}
(\cf~\cite{Kapranov}, (1.3)). In terms of free actions of roots of unity on
sets, and for $m=na$, the functor of restriction of scalars for the extension
\eqref{fieldext} is  obtained by viewing $\Z/n\Z$ as the subgroup of $ \Z/m\Z$
generated by $a$, which is in agreement with \eqref{morphism}.

Note also that there is no analogue for $q=1$ of the classification of finite
extensions of $\F_q$ for $q$ a prime power, and it is unjustified to consider
the inductive limit $\F_{1^\infty}=\varinjlim \F_{1^n}$  of the extensions
\eqref{fieldext} as the algebraic closure of $\F_1$.

\subsection{Affine varieties over $\F_1$}

We start by recalling briefly the notion of an affine variety over $\F_1$ as
introduced in \cite{Soule}.  Starting with the category of
 (commutative) rings $R$ with unit, which are finite and flat over $\Z$, we denote
 by  $\cR$ the full sub-category generated
 by the rings $A_n=\Z[T]/(T^n-1)$ as in \eqref{basic} and their tensor products
 (as $\Z$-modules). A gadget\footnote{``truc" in French} $X= (\underline X,\cA_X,e_X)$
 over $\F_1$ is specified (\cf~\cite{Soule}, Definition~3, \S~3.4) by
 giving the following data:\smallskip

 (a)~A covariant functor $\underline X: \cR\to \mathcal Sets$ to the category of  sets.

 (b)~A  $\C$-algebra $\cA_X$.

 (c)~A natural transformation $e_X$ from the functor $\underline X$
 to the functor $R\mapsto \Hom(\cA_X,R_\C)$.

 The notion of morphism of gadgets is that of natural transformation \ie a
morphism $\phi$ from $X$ to $Y$ is given by a pair $\phi=(\underline \phi,
\phi^*)$
\begin{equation}\label{fungadg}
    \underline \phi: \underline X\to \underline Y\,, \ \ \phi^*: \cA_Y\to \cA_X
\end{equation}
where $\underline \phi$ is a natural transformation of functors and $\phi^*$ a
morphism of algebras. One requires the compatibility with the evaluation maps
\ie one has a commutative diagram
\begin{gather}
\label{functmap}
 \,\hspace{50pt}\raisetag{-47pt} \xymatrix@C=25pt@R=25pt{
 \underline X(R)\ar[d]_{e_X(R)}\ar[r]^-{\underline \phi(R)} &
  \underline Y(R)\ar[d]^{e_Y(R)}& \\
\Hom(\cA_X,R_\C) \ar[r]^-{\phi^*}  & \Hom(\cA_Y,R_\C)\\
}\nonumber
\end{gather}
A morphism $\phi$ from $X$ to $Y$ is an {\em immersion} when $\phi^*$ is
injective and for any object $R$ of $\cR$, the map $ \underline \phi(R):
\underline X(R)\to \underline Y(R) $ is injective.

The construction described in \S 3.3 of
  \cite{Soule}  gives a natural covariant functor $\truc$ from the category of
  varieties over $\Z$ (\ie schemes of finite type over $\Z$) to
  the category $\cT$ of gadgets over $\F_1$. More precisely

  \begin{lem}\label{functorF}
  An affine variety $V$ over ${\mathbb Z}$ defines a gadget $X=\cF(V)$ over
${\mathbb F}_1$ by letting
 \begin{equation}\label{trucV}
    \underline X(R)={\rm Hom}_\Z(\cO(V),R)\,, \ \ \cA_X=\cO(V)\otimes_\Z \C, \ \ e_X(f)=f\otimes
    id_\C, \forall f\in \underline X(R)\,.
  \end{equation}
  \end{lem}

  One then defines (\cf \cite{Soule} Definition 3)

  \begin{defn}\label{varoverfundefn}
An  {\em affine} variety over $\F_1$ is a gadget $X$ such that:

 $\bullet$ $X(R)$
is finite for any object $R$ of $\cR$.

$\bullet$ There exists an affine variety $X_{\mathbb Z}$ over ${\mathbb Z}$ and
an immersion $i : X \rightarrow \cF(X_{\mathbb Z})$ of gadgets satisfying the
following property: for any affine variety $V$ over ${\mathbb Z}$ and any
morphism of gadgets $ \varphi : X \rightarrow \cF(V)\, , $ there exists a
unique algebraic morphism
$$
\varphi_{\mathbb Z} : X_{\mathbb Z} \rightarrow\textbf{} V
$$
such that $\varphi = \cF(\varphi_{\mathbb Z}) \circ i$.
\end{defn}

\medskip

\subsection{The varieties $\mu^{(k)}$}

We introduce the varieties $\mu^{(k)}$ over $\F_1$, as examples of  affine
varieties over $\F_1$. We begin by defining the associated functors
${\underline \mu}^{(k)}: \cR \to \mathcal
   Sets$. These are given (for $k\in \N$) by the assignments
\begin{equation}\label{deffunmuk}
{\underline \mu}^{(k)}(R)=\{x\in R\,|\,x^k=1\},
\end{equation}
\ie  ${\underline \mu}^{(k)}$ is the functor  represented by the ring $A_k$
\begin{equation}\label{funmuk2}
    {\underline \mu}^{(k)}(R)={\rm Hom}_\Z(A_k,R)\qqq R\in \mathcal Obj(\cR)\,.
\end{equation}
Notice that the functors ${\underline \mu}^{(k)}$ can be organized in two
different ways:
\begin{enumerate}
\item[a)] As an inductive system converging to ${\underline \mu} = \underline{\Gm}$,
where ${\Gm}(R):=\mu(R)$ is the multiplicative group $\Gm$ over $\F_1$, as affine variety over $\F_1$ (\cf~\cite{Soule}, 5.2.2).
\item[b)] As a projective system converging to ${\underline \mu}^{(\infty)}$, where \begin{center} ${\underline \mu}^{(\infty)}(R):= \varprojlim_k{\rm Hom}(A_k,R) = {\rm
  Hom}(\Z[\Q/\Z],R)$.\end{center}
\end{enumerate}

For a), one uses the natural inclusion
\begin{equation}\label{incl}
    {\underline \mu}^{(n)}(R)\subset {\underline \mu}^{(m)}(R)\quad \forall~n|m
\end{equation}
which corresponds at the level of the rings $A_n$ representing these functors,
to the surjective ring homomorphism
\begin{equation}\label{quorel}
A_m\twoheadrightarrow A_n=A_m/(T^n-1)\quad \forall~n|m.
\end{equation}
Then, the union of the ${\underline \mu}^{(k)}(R)$ is simply the functor
${\underline \mu}(R)$ which assigns to $R\in \mathcal Obj(\cR)$ the set of all
roots of $1$ in $R$. In the formalism of \cite{Soule} this functor is part (a)
of the data (a)-(c) for the multiplicative group $\mu={\Gm}$ as an affine
variety over $\F_1$.

\smallskip
To explain  b), we use the homomorphisms \eqref{morphism}. These homomorphisms
organize the ${\underline \mu}^{(k)}(R)$ into a projective system. The
connecting maps are given by raising a root of $1$ to the power $a=m/n$. Then,
the elements of ${\underline \mu}^{(\infty)}(R)$ are described by homomorphisms
of the group   $\Q/\Z$ to the multiplicative group of $R$. The equality
${\underline \mu}^{(\infty)}(R)={\rm
  Hom}(\Z[\Q/\Z],R)$ follows from \eqref{funmuk2}.
\smallskip

After tensoring by $\Z$ as in \eqref{basic}, the scalars extensions
$\F_{1^n}\subset \F_{1^m}$ of \eqref{fieldext} (\cf~\cite{Kapranov}, (1.3))
correspond to homomorphisms of rings
\begin{equation}\label{fieldext1}
  \xi_{m,n}:  \F_{1^n}\otimes_{\F_1}\Z\to \F_{1^m}\otimes_{\F_1}\Z\,
\end{equation}
given by
\begin{equation}\label{morphismF1}
    \xi_{m,n}(u(n))=u(m)^a\,, \ \ a=m/n\,,
\end{equation}
where $u(n)$ is the canonical generator $T\in A_n$. These agree with the maps
\eqref{morphism} that define the integral version of the abelian part of the
BC-system.

\medskip

In order to complete the definition of the varieties $\mu^{(n)}$ over $\F_1$,
we use  the functor $\cF$ of Lemma \ref{functorF}. In other words we
  {\em define}
  \begin{equation}\label{defnmun}
    \mu^{(n)}=\truc \Sp(A_n).
  \end{equation}
 One checks (\cf \cite{Soule}, \S 4, Proposition 2)
 that it fulfills the conditions of Definition \ref{varoverfundefn}.
 We then obtain the following result.

\begin{prop}\label{mukF1}
The $\mu^{(n)}$ defined as in \eqref{defnmun} form a projective system of
zero-dimensional affine varieties over $\F_1$.
\end{prop}

\proof
 It follows from \eqref{funmuk2} that the corresponding functors ${\underline
  \mu}^{(k)}$
  are the same as the ones defined by \eqref{deffunmuk}.
  The morphisms \eqref{morphism} turn the
  varieties $\Sp(A_n)$ into a projective system and thus, since
  $\truc$ is a covariant functor, we get that the $\mu^{(n)}$ form a
  projective system of varieties over $\F_1$.
\endproof

\section{The integral BC-endomotive}\label{BCintegral}

Having to work over $\Z$ creates a problem when one implements the semigroup
action via the maps $\rho_n$, which involve denominators, as in \eqref{endo}
and \eqref{idem}. However, as shown in the algebro-geometric description of the
BC-algebra (Proposition~\ref{propendo}, (d)), the partial inverses of the
$\rho_n$, which we have denoted by $\sigma_n$, do not involve denominators,
therefore we will be able to consider them over $\Z$.

The partial inverse relations between the $\sigma_n$ and $\rho_n$ are given by
\eqref{rhosigma1} and \eqref{rhosigma2}.

Since by \eqref{funmuk2} the schemes $\underline{\mu}^{(n)}$ are represented by
the rings $A_n$, by Yoneda's lemma the ring homomorphisms $\sigma_n$ given by
 \begin{equation}\label{actionlevelk}
    \sigma_n: A_k \to A_k\, \ \ \ \ u(k)\mapsto u(k)^n
  \end{equation}
define (contravariantly) morphisms of schemes. These induce morphisms of the
pro-scheme $\underline{\mu}^{(\infty)}$ by the compatibility
\begin{equation}\label{compproj}
    \xi_{\ell,k}\circ \sigma_n=\sigma_n\circ \xi_{\ell,k} \qqq
    k|\ell \qqq n\in \N.
\end{equation}
In Proposition \ref{propdeffun} we will show how the maps $\sigma_n$ in fact
give rise to endomorphisms of the varieties $\mu^{(k)}$ over $\F_1$.

 In the limit, the endomorphisms $\sigma_n$ are surjective
\begin{equation}\label{surjectssigma}
    \sigma_n\,:\, \varinjlim_k A_k \twoheadrightarrow \varinjlim_k A_k .
\end{equation}
In fact, in the group ring notation of \S \ref{BCrecall1}, one gets
$\sigma_n(e(r))=e(nr)$, while one has the surjectivity of multiplication by $n$
in the exact sequence \eqref{absequ}. However, the $\sigma_n$ are {\em not} the
same as the endomorphisms $\rho_n$, since the latter are injective and describe
ring isomorphisms between reduced algebras and the projectors $\pi_n$, as we
have shown in Proposition~\ref{propendo}.

The kernel of $\sigma_n$ in \eqref{surjectssigma} is the ideal $J_n$ generated
by the element $u(n)-1$, or in group-ring notation by $e(1/n)-1$. This means
that $\sigma_n^{-1}(f)$ is only defined modulo $J_n$. If one allows inverting
$n$, then there is a natural complementary subspace to $J_n$, namely the
reduced algebra by the projection $\pi_n$. However, when we work over $\Z$ (and
a fortiori over $\F_1$) we cannot invert $n$, and we need to adapt the
presentation of the BC-endomotive. The data of the BC-endomotive, \ie the
abelian algebra and the endomorphisms, combine to produce a noncommutative
crossed product algebra with a natural time evolution defined over $\C$. This
quantum statistical mechanical system is the BC-system which we recall below.

\smallskip

\subsection{$C^*$-algebra description of the BC-system}\label{BCcstar}
 At the $C^*$-algebra level the BC system is given by
 $1_{\hat \Z}(C_0(\A_{\Q,f})\rtimes \Q_+^*)1_{\hat \Z}$, namely by the
 algebra of the crossed product $C_0(\A_{\Q,f})\rtimes \Q_+^*$ reduced by the
projection $1_{\hat \Z}\in C_0(\A_{\Q,f})$. Here $\A_{\Q,f}$ denotes the
locally compact space of finite adeles of $\Q$ and $\hat \Z\subset \A_{\Q,f}$
the open compact subset closure of $\Z$. The reduced algebra can be described
as the convolution algebra of the locally compact \'etale groupoid $\cG$
obtained as the reduction of the groupoid $\A_{\Q,f}\rtimes \Q_+^*$ by the open
and closed set of units $\hat \Z\subset \A_{\Q,f}$. Concretely, the groupoid
$\cG$ is the \'etale groupoid of pairs
\begin{equation}\label{G1pairs}
\cG=\{ (r,\rho)\, | \, r\in \Q^*_+,\, \rho\in \hat\Z, \, \text{ such that }
r\rho\in \hat\Z \},
\end{equation}
with source and range maps $(r,\rho)\mapsto \rho$ and $(r,\rho)\mapsto r\rho$,
and composition
\begin{equation}\label{cG1comp}
(r_1,\rho_1)\circ (r_2,\rho_2)=(r_1 r_2,\rho_2), \ \ \ \text{ if }
r_2\rho_2=\rho_1.
\end{equation}
 The $C^*$-algebra $C^*(\cG)$ of a locally compact \'etale groupoid
$\cG$ is obtained as the completion of the algebra $C_c(\cG)$ of compactly
supported functions on $\cG$ with the convolution product
\begin{equation}\label{cGconvol}
f_1*f_2(g)=\sum_{g_1g_2=g} f_1(g_1)f_2(g_2),
\end{equation}
 the involution
\begin{equation}\label{involution1}
f^*(g)=\overline{f(g^{-1})} \end{equation}
 and the norm
\begin{equation}\label{cGnorm}
\| f \|:= \sup_{y\in \cG^{(0)}} \| \pi_y (f)\|_{\cB(\cH_y)}.
\end{equation}
Here every unit $y\in \cG^{(0)}$ defines a representation $\pi_y$ by left
convolution of the algebra $C_c(\cG)$  on the Hilbert space
$\cH_y=\ell^2(\cG_y)$, where $\cG_y$ denotes the set of elements in $\cG$ with
source $y$. Namely, one has
\begin{equation}\label{cGpiy}
(\pi_y(f)\xi)(g) = \sum_{g_1g_2=g} f(g_1)\xi(g_2).
\end{equation}
 The $C^*$-algebra $C^*(\cG)$ contains $C(\hat \Z)$ as a subalgebra and is
generated by $C(\hat \Z)$ and the elements $\mu_n$  given by the compactly
supported functions
\begin{equation}\label{muen}
    \mu_n(n,\rho)=1\qqq \rho\in \hat\Z\,,  \ \ \mu_n(r,\rho)=0 \qqq r\neq n\,, \rho\in
    \hat\Z\,.
\end{equation}
One identifies
 the Pontrjagin dual of the group $\Q/\Z$ with the compact group $\hat\Z= \Hom(\Q/\Z,\Q/\Z)$ using the
 pairing $$\langle
    \gamma,\rho\rangle= e^{2\pi i \rho(\gamma)}\qqq \gamma \in \Q/\Z\,, \ \rho \in \Hom(\Q/\Z,\Q/\Z) $$ and one
 lets $e(\gamma)\in C^*(\cG)$ be given by the function
 \begin{equation}\label{egamma}
    e(\gamma)(r,\rho)=0 \qqq r\neq 1\,,\ \ e(\gamma)(1,\rho)=\langle
    \gamma,\rho\rangle\,.
 \end{equation}

The  time evolution is given by the following one-parameter group of
automorphisms of the $C^*$-algebra $C^*(\cG)$:
\begin{equation}\label{sigmatgenerators}
\sigma_t (\mu_n) = n^{it} \mu_n , \ \ \sigma_t (\mu_n^*) = n^{-it} \mu_n^* , \
\ \ \ \ \sigma_t (e(\gamma))=e(\gamma).
\end{equation}\vspace{.001in}

\begin{defn} \label{bcsystem} The BC-system is the complex dynamical system defined by the pair
$(C^*(\cG),\sigma_t)$.
\end{defn}

 We refer to \cite{CMbook} Chapter 3, \S 4 for the equivalent descriptions of
 the $C^*$-algebra of the BC-system and of the relation with $\Q$-lattices. Working over $\C$ one  considers
 the subalgebra of $C^*(\cG)$ generated by the characters $e(\gamma)$, $\gamma\in \Q/\Z$ the $\mu_n$ and their
 adjoints $\mu_n^*$. We
 shall now explain the presentation of this algebra over $\Q$.
\smallskip

\subsection{The BC-algebra over $\Q$}\label{bcqsect}

We first recall the presentation of the crossed product algebra
$\cA_\Q=\Q[\Q/\Z]\rtimes \N$ of the BC system in characteristic zero.

The group ring $\Q[\Q/\Z]$ has the canonical additive basis  $e(\gamma)$,
$\gamma\in\Q/\Z$, with $e(\gamma)^* = e(-\gamma)$ and $e(\gamma_1+\gamma_2) =
e(\gamma_1)e(\gamma_2)$. To obtain the crossed product, one considers then
generators $\mu_n$ and $\mu_n^*$, $n\in\N$, which satisfy the following
conditions:\medskip

(c1)~$\mu_n^*\mu_n = 1\,, \ \ \forall~n$\smallskip

(c2)~$\mu_{nm} = \mu_n\mu_m\,, \ \mu^*_{nm} = \mu^*_n\mu^*_m\,, \  \ \forall~n,m$\,,\smallskip

(c3)~$\mu_n\mu_m^* = \mu_m^*\mu_n\,, \ \ \text{if}~ (m,n)=1$\medskip

together with the additional relation\medskip

(c4)~$\mu_n e(\gamma)\mu_n^* = \frac{1}{n}\sum_{n\delta =
\gamma}e(\delta)\,, \ \ \forall n,~\gamma$.\medskip

In particular, the relation (c4) can be interpreted algebraically by means of
the homomorphism $\rho_n(x)$ (\cf~\eqref{endo}) projecting onto the reduced
algebra by the idempotent $\pi_n$. This means\medskip

(c4')~$\rho_n(x) = \mu_n x\mu_n^*\,, \ \ \forall~x\in\Q[\Q/\Z]$.\medskip

As a consequence of (c1) and (c4') we get
\[
(\mu_n\mu_n^*)^2=(\mu_n\mu_n^*)(\mu_n\mu_n^*) = \mu_n\mu_n^* = \rho_n(1) =
\frac{1}{n}\sum_{n\gamma=0}e(\gamma)
\]
In this way we get a description of the projector $\pi_n\in\Q[\Q/\Z]$ as in
\eqref{idem} by means of the new generators of the crossed product \ie
$\pi_n=\rho_n(1)=\mu_n\mu_n^*$. It also follows from (c1) that
$\mu_n\mu_n^*\mu_n = \mu_n$. Since the surjective endomorphisms $\sigma_n$ are
partial inverses of $\rho_n$, that is
\[
\sigma_n\rho_n(e(\gamma)) =
\frac{1}{n}\sum_{n\gamma'=\gamma}\sigma_n(e(\gamma')) =
\frac{1}{n}\sum_{n\gamma'=\gamma}e(n\gamma') = e(\gamma),
\]
one gets $\sigma_n\rho_n(x) = x$, $\forall~x\in\Q[\Q/\Z]$.\medskip

We then have the following easy consequence

\begin{prop}\label{murhosigma}
The following relations hold in the algebra $\cA_\Q=\Q[\Q/\Z]\rtimes_\rho\N$:
\begin{equation}\label{relmurho}
\mu_n x = \rho_n(x) \mu_n, \ \ \  \ \forall x \in \Q[\Q/\Z], \, \forall n \in
\N,
\end{equation}
\begin{equation}\label{relmusigma}
\mu_n^* x = \sigma_n(x) \mu_n^*, \ \ \  \ \forall x \in \Q[\Q/\Z], \, \forall n
\in \N,
\end{equation}
\begin{equation}\label{relmusigma2}
x \mu_n = \mu_n \sigma_n(x),  \ \ \  \ \forall x \in \Q[\Q/\Z], \, \forall n
\in \N.
\end{equation}
\end{prop}

\proof Relation \eqref{relmurho} follows from
\begin{equation}\label{mumurho}
\mu_n x \mu_n^* = \rho_n(x)
\end{equation}
and the fact that $\mu_n^* \mu_n =1$. For \eqref{relmusigma}, we use the
idempotent $\pi_n=\mu_n\mu_n^*=\rho_n(1)$. We first assume that $x =\pi_n x$
belongs to the reduced algebra by $\pi_n$. It then follows that $x= \rho_n(y)$
for some $y\in \Q[\Q/\Z]$ (hence $y=\sigma_n(x)$). By applying (c4') and (c1),
this shows that
\[
\mu_n^* x= \mu_n^* \rho_n(y) = \mu_n^* \mu_n y \mu_n^*= y \mu_n^* =\sigma_n(x)
\mu_n^*.
\]
In the general case, we notice that in view on (c1), the left hand side of
\eqref{relmusigma} does not change by replacing $x$ by $\mu_n\mu_n^* x= \pi_n
x$. The right hand side does not change either, since $\sigma_n (\pi_n) =1$,
hence  \eqref{relmusigma} holds with no restriction. The relation
\eqref{relmusigma} also gives
\begin{equation}\label{mustarmusigma}
\mu_n^*x\mu_n = \sigma_n(x)
\end{equation}
by multiplying on the right by $\mu_n$ and applying (c1). The relation
\eqref{relmusigma2} then follows by \eqref{mustarmusigma} together with $x\mu_n
= \pi_n x \mu_n = x \pi_n \mu_n$.
\endproof

\begin{rem}{\rm Notice that the involution \eqref{involution1} of the
$C^*$-algebra $C^*(\cG)$ restricts to an involution of the rational algebra
$\cA_\Q=\Q[\Q/\Z]\rtimes \N$ with the properties
\begin{equation}\label{involrat}
    e(\gamma)\mapsto e(-\gamma)\,, \ \ \mu_n\mapsto \mu_n^*\,, \ \
    \mu_n^*\mapsto \mu_n\,.
\end{equation}
Note that the full presentation of the rational algebra involves the two
relations that appear in (c2). In particular this is needed for the involution
\eqref{involrat} to make sense.  }\end{rem}

\medskip

\subsection{The maps $\tilde\rho_n$}\label{secttilderho}

When one wants to generalize the definition of the algebra
$\cA_\Q=\Q[\Q/\Z]\rtimes_\rho\N$ to the case where the field of coefficients is
a perfect field $\K$ of positive characteristic (for example $\K = \F_p$), as
well as in extending the original (rational) formulation of the algebra to the
case of integer coefficients, one is faced with the problem of ``dividing by
$n$" in the definition of the endomorphisms $\rho_n$ (\eg when $\K=\F_p$, for
$n=p$). However,  up to multiplying  the original definition of the maps
$\rho_n$ by $n$, the linear maps
\begin{equation}\label{newrhos}
\tilde\rho_n: \K[\Q/\Z] \to \K[\Q/\Z]\,, \ \ \tilde\rho_n(e(\gamma))=
\sum_{n\gamma'=\gamma}e(\gamma')
\end{equation}
retain a meaning (when $\text{char}(\K)=p>0$ and $n=p$ and also over $\Z$),
since $e(\gamma)$ is a linear basis of the algebra $\K[\Q/\Z]$ as a $\K$-vector
space. In particular, when $\text{char}(\K)=p>0$, the operator $\tilde{\pi}_p =
\tilde\rho_p(1) = \sum_{p\gamma=0}e(\gamma) \neq 0$ is nilpotent since
\[
\tilde{\pi}_p^2 =
(\sum_{p\gamma=0}e(\gamma))(\sum_{p\gamma'=0}e(\gamma'))=\sum_{\gamma,\gamma'}e(\gamma+\gamma')
= p\sum_{p\gamma''=0}e(\gamma'') = 0.
\]
Compare this with the idempotents $\pi_n$ of \eqref{idem}.  Moreover, over a
perfect field of characteristic $p>0$ one gets $\sigma_p\tilde\rho_p = 0$,
since $\sigma_p\tilde{\pi}_p = \sigma_p\tilde\rho_p(1) =
\sigma_p\sum_{p\gamma=0}e(\gamma) = \sum_{p\gamma=0}e(p\gamma) = 0$. This means
that $\text{Range}(\tilde\rho_p)\subset\text{Ker}(\sigma_p)$.

\begin{prop}\label{mapstilderho} When working over $\Z$ the $\sigma_n$ continue to make sense and define
endomorphisms of $\Z[\Q/\Z]$ which fulfill the following relations with the
maps $\tilde\rho_m$:
\begin{equation}\label{rhomult}
 \sigma_{nm}=\sigma_{n}\sigma_{m}\,,\ \
    \tilde\rho_{mn}=\tilde\rho_{m}\tilde\rho_{n}\qqq m,n
\end{equation}
\begin{equation}\label{rhosigmamult}
    \tilde\rho_{m}(\sigma_m(x)y)=x\tilde\rho_{m}(y)\qqq x,y\in \Z[\Q/\Z]
\end{equation}
\begin{equation}\label{divisible}
    \sigma_c(\tilde\rho_b(x))=(b,c)\,\tilde\rho_{b'}(\sigma_{c'}(x))\,, \ \
    b'=b/(b,c)\,, \ \ c'=c/(b,c)\,,
\end{equation}
where $(b,c)$ denotes the gcd of $b$ and $c$.
\end{prop}

\proof One has by definition $\sigma_n(e(\gamma))=e(n\gamma)$, which shows that
$\sigma_n$ is an endomorphism of $\Z[\Q/\Z]$ and
$\sigma_{nm}=\sigma_{n}\sigma_{m}$. To get
$\tilde\rho_{mn}=\tilde\rho_{m}\tilde\rho_{n}$ we let, for $x\in \Q/\Z$ and
$n\in \N$,
\begin{equation}\label{En}
E_n(x)=\{y\in \Q/\Z\,|\, ny=x\}\,.
\end{equation}
 One has
$$
E_{nm}(x)=\cup_{y\in E_n(x)}E_m(y)\,, \ \ y_1\neq y_2\Rightarrow E_m(y_1)\cap
E_m(y_2)=\emptyset
$$
thus
$$
\tilde\rho_{m}(\tilde\rho_{n}(e(x)))=\tilde\rho_{m}(\sum_{E_n(x)}e(y))=\sum_{y\in
E_n(x)}\sum_{z\in E_m(y)}e(z)=\tilde\rho_{mn}(e(x))\,.
$$
To check \eqref{rhosigmamult} we can assume that $x=e(s)$, $y=e(t)$ with
$s,t\in \Q/\Z$. One has $\sigma_m(x)y=e(ms)e(t)=e(ms+t)$. For $u\in \Q/\Z$, one
has $mu=ms+t$ iff $u-s\in E_m(t)$ thus $E_m(ms+t)=s+E_m(t)$ which proves
\eqref{rhosigmamult}.

To check \eqref{divisible} we  assume that $x=e(s)$ and let $n=(b,c)$ so that
$b=nb'$, $c=nc'$ with $(b',c')=1$. One has
$$
E_b(s)=\{u\in \Q/\Z\,|\, bu=s\}=\{u\in \Q/\Z\,|\, nb'u=s\}
$$
Thus the multiplication by $c=nc'$ is an $n$ to $1$ map from $E_b(s)$ to
$E_{b'}(c's)$. This proves \eqref{divisible}.
\endproof

In particular one gets:

\begin{cor}\label{mapstilderhocor} The range of $\tilde\rho_{m}$ is an ideal in
$\Z[\Q/\Z]$. When $n$ and $m$ are relatively prime $\sigma_n$ commutes with
$\tilde\rho_{m}$.
\end{cor}

\proof The range of $\tilde\rho_{m}$ is additive by construction and is
invariant under multiplication by $\Z[\Q/\Z]$ using \eqref{rhosigmamult}. The
second statement follows from \eqref{divisible}.
\endproof

\begin{rem}{\rm Notice that, although the $\tilde \rho_n$ are not ring
homomorphisms, the relation \eqref{rhosigmamult} which they fulfill suggests
the existence of an associated correspondence (in the form of a bimodule). This
would fit with a more general framework for the theory of endomotives that uses
correspondences instead of endomorphisms as in \cite{Pimsner}. }\end{rem}

\subsection{The BC-algebra over $\Z$}\label{bczz}
When $\K$ denotes either $\Z$ or a perfect field of positive characteristic,
the relations  \eqref{relmusigma} continue to make sense, because the
$\sigma_n$ are well defined. On the other hand, the relation \eqref{relmurho}
involves the $\rho_n$ which are not well defined. However, in the case of
integral coefficients and in characteristic $p$,
 the linear maps $\tilde\rho_n$ of \eqref{newrhos}
make sense and in the latter case these maps play the role of the $p\rho_p$.
Thus, in order to extend the relation \eqref{relmurho},  we keep the generators
$\mu_n^*$ and introduce new generators $\tilde\mu_n$ (in place of the $\mu_n$'s), which play the role, in characteristic $p$, of
the operators $p\mu_p$ and  in general fulfill the relation
\begin{equation}\label{pmuprel}
\tilde\mu_n x \mu_n^* = \tilde\rho_n(x),
\end{equation}
that is the analog of \eqref{mumurho}. These relations reformulate (c4') in the
case of integral coefficients and make sense in positive characteristic.

\begin{defn} \label{overZ} The algebra
$\cA_\Z=\Z[\Q/\Z]\rtimes_{\tilde\rho}\N$ is the algebra generated by the group
ring $\Z[\Q/\Z]$, and by the elements $\tilde\mu_n$ and $\mu_n^*$, with
$n\in\N$,
 which satisfy
the relations:
\begin{equation}\label{presoverZ1}
\begin{array}{l}
\tilde\mu_n x \mu_n^* = \tilde\rho_n(x)\ \ \ \  \\[3mm]
\mu_n^* x = \sigma_n(x) \mu_n^*  \\[3mm]
x \tilde\mu_n = \tilde\mu_n \sigma_n(x),
\end{array}
\end{equation}
where $\tilde\rho_m$, $m\in\N$ is defined in \eqref{newrhos}, as well as the
relations
\begin{equation}\label{presoverZ2}
\begin{array}{l}
\tilde\mu_{nm}= \tilde\mu_n  \tilde\mu_m \qqq n,m\\[3mm]
\mu_{nm}^* =\mu_{n}^*\mu_{m}^* \qqq n,m\\[3mm]
\mu_n^* \tilde\mu_n  =n \\[3mm]
\tilde\mu_n\mu_m^* = \mu_m^*\tilde\mu_n \ \ \ \ (n,m) = 1.
\end{array}
\end{equation}
\end{defn}

Our first task is to check that these relations are sufficient to express every
element of $\cA_\Z$ as a finite sum of elementary monomials labeled by a pair
$(x,r)$ where $x\in \Z[\Q/\Z]$ and $r\in \Q_+^*$ is an irreducible
fraction $r=a/b$.

\begin{lem}\label{pres1Z} Any element of  the
algebra $\cA_\Z$ is a finite sum of monomials:
\begin{equation}\label{basismonZ}
    \tilde\mu_a\,x \, \mu^*_b\,, \ \ \ (a,b) = 1\,, \  \ x\in \Z[\Q/\Z]\,,
\end{equation}
where by convention $\tilde\mu_1=\mu^*_1=1$.
\end{lem}

\proof It is enough to show that the product of monomials of the form
\eqref{basismonZ} is still of the same form. We take a product of the form
$$
\tilde\mu_a\,x \, \mu^*_b\,\tilde\mu_c\,y \, \mu^*_d\,
$$
Let then $n$ be the gcd of $b=nb'$ and $c=nc'$. One has
$$
\mu^*_b\,\tilde\mu_c=\mu^*_{b'}\,\mu^*_n\,\tilde\mu_n\tilde\mu_{c'}=n
\mu^*_{b'}\,\tilde\mu_{c'}=n\,\tilde\mu_{c'}\mu^*_{b'}
$$
so that
$$
\tilde\mu_a\,x \, \mu^*_b\,\tilde\mu_c\,y \, \mu^*_d=n\,\tilde\mu_a\,x \,
\tilde\mu_{c'}\mu^*_{b'}\,y \,
\mu^*_d=\,n\,\tilde\mu_a\,\tilde\mu_{c'}\sigma_{c'}(x)\sigma_{b'}(y)
\mu^*_{b'}\, \mu^*_d
$$
Let then $m$ be the gcd of $ac'=mu$ and $b'd=mv$. One has
$$
\tilde\mu_a\,\tilde\mu_{c'}=\tilde\mu_u\,\tilde\mu_{m}\,, \ \ \mu^*_{b'}\,
\mu^*_d=\mu^*_{m}\, \mu^*_v
$$
so that:
$$
\tilde\mu_a\,x \, \mu^*_b\,\tilde\mu_c\,y \,
\mu^*_d=n\,\tilde\mu_u\,\tilde\mu_{m}\,\sigma_{c'}(x)\sigma_{b'}(y)\,\mu^*_{m}\,
\mu^*_v= n\,\tilde\mu_u\,
\tilde\rho_m(\sigma_{c'}(x)\sigma_{b'}(y))\,\mu^*_v\,.
$$
Since $u$ and $v$ are relatively prime and
$z=\tilde\rho_m(\sigma_{c'}(x)\sigma_{b'}(y))\in \Z[\Q/\Z]$ it follows that the
product of two monomials of the form \eqref{basismonZ} is still a monomial of
the same form. Note also that
\begin{equation}\label{multrat}
    u/v=(a/b)(c/d)
\end{equation}
since $u/v=ac'/(b'd)=ac/(bd)$. Thus the labels $a/b\in \Q_+^*$ are
multiplicative.
\endproof

\begin{rem}\label{geomreason} {\rm Using the surjectivity of the endomorphisms
$\sigma_n$ one can rewrite the monomials \eqref{basismonZ} in the form
$y\,\tilde\mu_a\, \mu^*_b$, $\tilde\mu_a\, \mu^*_b\, z$ or $\mu^*_b\, t\,
\tilde\mu_a$. The reason for choosing \eqref{basismonZ} is that, in this form,
there is no ambiguity in the choice of $x$ while the lack of injectivity of
$\sigma_a$ and $\sigma_b$ introduces  an ambiguity in the choices of $y$, $z$
and $t$. At the geometric level this corresponds, using \eqref{muen}, to the
fact that the initial support of $\tilde\mu_a$ is $1$. } \end{rem}

In order to check that the relations of Definition \ref{overZ} are coherent we
shall now construct a faithful representation of these relations (which is the
left regular representation of $\cA_\Z$) in the free abelian group
$\cE=\Z[\Q/\Z\times \Q_+^*]$. We denote by $\xi(x,r)$ the element of $\cE$
associated to $x\in \Z[\Q/\Z]$ and $r\in  \Q_+^*$.

\begin{prop} \label{pres4} The following relations
define a faithful  representation of the algebra $\cA_\Z$ on $\cE$,
\begin{equation}\label{leftactZ1}
    x \,\xi(y,c/d)=\xi(\sigma_c(x)y,c/d)\qqq c,d\,, \ (c,d)=1
\end{equation}
\begin{equation}\label{leftactZ2}
\tilde\mu_a\,\xi(y,c/d)= \xi(\tilde\rho_m(y),ac/d)\,, \ \ m=(a,d)
\end{equation}
\begin{equation}\label{leftactZ3}
\mu_b^*\,\xi(y,c/d)=(b,c)\, \xi(\sigma_{b/n}(y),c/bd)\,, \ \ n=(b,c)\,.
\end{equation}
\end{prop}

\proof We shall check that  the relations of Definition \ref{overZ} are
fulfilled. The relation \eqref{leftactZ1} shows that the left action of
$\Z[\Q/\Z]$ is a representation which is a direct sum of copies of the left
regular representation of $\Z[\Q/\Z]$ composed with the $\sigma_c$.

Using the notation $(a,b)$ for ${\rm gcd}(a,b)$ one has the equality
\begin{equation}\label{gcd}
    (a_1a_2,d)=(a_1,d)(a_2,d/(a_1,d))
\end{equation}
and the fact that the left action of $\tilde\mu_a$ fulfills
$\tilde\mu_{a_2a_1}=\tilde\mu_{a_2}\tilde\mu_{a_1}$ follows from
\eqref{rhomult} which gives
$$
\tilde\rho_m(y)=\tilde\rho_{m_2}(\tilde\rho_{m_1}(y))\,, \ m_1=(a_1,d)\,, \
m_2=(a_2,d/(a_1,d))\,, \ m=(a_1a_2,d)\,.
$$
In order to check the relation $x \tilde\mu_a = \tilde\mu_a \sigma_a(x)$ one
uses \eqref{rhosigmamult}. One has
$$
x \tilde\mu_a\,\xi(y,c/d)=x
\,\xi(\tilde\rho_m(y),ac/d)=\xi(\sigma_k(x)\tilde\rho_m(y),ac/d)\,, \ k=ac/m\,,
\ m=(a,d)
$$
$$
 \tilde\mu_a \sigma_a(x)\,\xi(y,c/d)=\tilde\mu_a
 \xi(\sigma_c(\sigma_a(x))y,c/d)= \xi(\tilde\rho_m(\sigma_{ac}(x)y),ac/d)
$$
and since $ac=mk$, \eqref{rhosigmamult} gives
$$
\tilde\rho_m(\sigma_{ac}(x)y)=\sigma_k(x)\tilde\rho_m(y)\,.
$$
Let us check  the relation $\mu^*_{b_2b_1}=\mu^*_{b_2}\mu^*_{b_1}$. Let
$n_1=(b_1,c)$ and $b'_1=b_1/n_1$, $c'_1=c/n_1$ then
$$
\mu_{b_1}^*\,\xi(y,c/d)=n_1\,\xi(\sigma_{b'_1}(y),c'_1/(b'_1d))
$$
so that, with $n_2=(b_2,c'_1)$ and $b'_2=b_2/n_2$, $c'_2=c'_1/n_2$ one gets
$$
\mu_{b_2}^*(\mu_{b_1}^*\,\xi(y,c/d))=n_2n_1\,
\xi(\sigma_{b'_2}\sigma_{b'_1}(y),c'_2/(b'_2b'_1d))\,.
$$
By \eqref{gcd} one has $n_1n_2=(b_1,c)(b_2,c/n_1)=(b_1b_2,c)=n$ and with
$b=b_1b_2$ one has $$b'=b/n=(b_1/n_1)(b_2/n_2)=b'_1b'_2$$
$$c'_2=c'_1/n_2=c/(n_1n_2)=c/n=c'$$
This shows, using $\sigma_{ab}=\sigma_a\sigma_b$, that
$$
\mu_{b_2}^*(\mu_{b_1}^*\,\xi(y,c/d))=\mu_{b_2b_1}^*\,\xi(y,c/d)\,.
$$
Let us now check the relation $\mu_b^* x = \sigma_b(x) \mu_b^* $. One has, with
$n=(b,c)$, $b'=b/n$, $c'=c/n$,
$$
\mu_b^*(
x\,\xi(y,c/d))=\mu_b^*\,\xi(\sigma_c(x)y,c/d)=n\,\xi(\sigma_{b'}(\sigma_c(x)y),c/bd)
$$
$$
\sigma_b(x)(\mu_b^*\,\xi(y,c/d))=n\,\sigma_b(x)\xi(\sigma_{b'}(y),c'/b'd)=n\,
\xi( \sigma_{c'}(\sigma_b(x))\,\sigma_{b'}(y),c/bd)
$$
Thus the relation follows from the multiplicativity of $\sigma_{b'}$ and the
equality $b'c=c'b$.

Let us check the relation $\tilde\mu_b x \mu_b^* = \tilde\rho_b(x)$. One has
$$
x (\mu_b^*(\xi(y,c/d)))=n\, \xi(\sigma_{c'}(x) \sigma_{b'}(y),c'/{b'd})\,, \ \
n=(b,c)\,, \ b'=b/n\,, \ c'=c/n\,.
$$
To multiply by $\tilde\mu_b$ on the left, one uses \eqref{leftactZ2} and gets
$$
\tilde\mu_b(x (\mu_b^*(\xi(y,c/d))))=n\,\xi( \tilde\rho_m(\sigma_{c'}(x)
\sigma_{b'}(y)),u/v)
$$
where $m=(b,b'd)$ and $u=bc'/m$, $v=b'd/m$. One has $m=b'$ since it divides
$b=n b'$ and $b'd$ while $bc'/b'=c$ is prime to $d$. Thus $u=c$ and $v=d$ and
one gets
$$
\tilde\mu_b(x (\mu_b^*(\xi(y,c/d))))=n\,\xi( \tilde\rho_{b'}(\sigma_{c'}(x)
\sigma_{b'}(y)),c/d)\,.
$$
In particular it is divisible by $n$ and one needs to understand why the other
side, namely $\tilde\rho_b(x)\xi(y,c/d)$ is also divisible by $n=(b,c)$. This
follows from \eqref{divisible} since,
$$
\tilde\rho_b(x)\xi(y,c/d)=\xi(\sigma_c(\tilde\rho_b(x))y,c/d)
$$
(by \eqref{leftactZ1}) while by \eqref{divisible},
$$
 \sigma_c(\tilde\rho_b(x))=(b,c)\,\tilde\rho_{b'}(\sigma_{c'}(x))\,, \ \
    b'=b/(b,c)\,, \ \ c'=c/(b,c)\,.
$$
 One then
uses \eqref{rhosigmamult} to obtain
$$
\tilde\rho_{b'}(\sigma_{c'}(x)
\sigma_{b'}(y))=\tilde\rho_{b'}(\sigma_{c'}(x))\,y
$$
which gives the required equality.

Let us now check the relation $\mu_a^*\tilde\mu_a=a$. By \eqref{leftactZ2} one
has
$$
\tilde\mu_a\,\xi(y,c/d)= \xi(\tilde\rho_m(y),u/v)
$$
where $m=(a,d)$ is the gcd of $ac=mu$ and $d=mv$. We then get with $a=ma'$,
$d=md'$ that $u=a'c$ and $v=d'$. The left action of $\mu_a^*$ is given by
$$
\mu_a^*(\tilde\mu_a\,\xi(y,c/d))=\mu_a^*\xi(\tilde\rho_m(y),u/v)=n\,\xi(
\sigma_{a''}(\tilde\rho_m(y)),c''/(a''d'))
$$
where $n=(a,a'c)$, $c''=a'c/n$, $a''=a/n$. One has $n=a'$ since $(m,c)=1$ as
$m=(a,d)$ is a divisor of $d$ and $(c,d)=1$. It follows that $c''=a'c/n=c$,
$a''=a/n=m$. Thus by \eqref{divisible} $\sigma_{a''}(\tilde\rho_m(y))=my$. Also
$a''d'=md'=d$, thus
$$
n\,\xi( \sigma_{a''}(\tilde\rho_m(y)),c''/(a''d'))=nm\,\xi(y,c/d)
$$
and the required equality follows from $nm=a$.

It remains to check that $\tilde\mu_a\mu_b^* = \mu_b^*\tilde\mu_a$ when $(a,b)
 = 1$. Let, as above,
 $m=(a,d)$ and write $a=ma'$,
$d=md'$ so that $u=a'c$ is prime to $v=d'$. One has
$$
\mu_b^*(\tilde\mu_a\,\xi(y,c/d))=\mu_b^*\,\xi(\tilde\rho_m(y),u/v)=n\,\xi(
\sigma_{b/n}(\tilde\rho_m(y)),ac/(bd))
$$
where $n=(b,u)$. Since $(a,b)=1$ one has $(a',b)=1$ and
$n=(b,u)=(b,a'c)=(b,c)$. Thus by \eqref{leftactZ3},
$$
\mu_b^*\,\xi(y,c/d)=n\,\xi(\sigma_{b/n}(y),c/bd)
$$
When applying $\tilde\mu_a$ on the left, one uses \eqref{leftactZ2}. One lets
$m=(a,b'd)$ where $b'=b/n$ so that $b'd$ is the reduced denominator of $c/bd$.
By \eqref{leftactZ2}, one has
$$
\tilde\mu_a\,\mu_b^*\,\xi(y,c/d)=n\,\tilde\mu_a\,\xi(\sigma_{b/n}(y),c/bd)=
n\,\xi(\tilde\rho_m(\sigma_{b/n}(y)),ac/bd)
$$
Since $(a,b)=1$ one has $m=(a,b'd)=(a,d)$ and the required equality follows
from the second statement of Corollary \ref{mapstilderhocor} since $m$ and
$b/n$ are relatively prime so that $\sigma_{b/n}$ and $\tilde\rho_m$ commute.
We have shown that the relations of Definition \ref{overZ} are fulfilled. One
has, for $(a,b)=1$,
\begin{equation}\label{toto}
    \tilde\mu_a\,x\,\mu_b^*\,\xi(1,1)=\xi(x,a/b)
\end{equation}
which shows that the map $x\in \cA_\Z\mapsto x\,\xi(1,1)\in \cE$ is an
isomorphism of abelian groups, and hence the representation of $\cA_\Z$ in
$\cE$  is faithful.
\endproof

\begin{cor} \label{corpres4}
The monomials
\begin{equation}\label{basismonZbis}
    \tilde\mu_a\,e(r) \, \mu^*_b\,, \ \ \ (a,b) = 1\,, \  \ r\in  \Q/\Z \,,
\end{equation}
form a basis of $\cA_\Z$ as a free abelian group.
\end{cor}

\proof By construction  $\cE=\Z[\Q/\Z\times \Q_+^*]$ is a free abelian group
with basis the $\xi(e(r),a/b)$ for $r\in  \Q/\Z$ and $a/b\in \Q_+^*$. Moreover
by \eqref{toto} and Proposition \ref{pres4} the map $x\in \cA_\Z\mapsto
x\,\xi(1,1)\in \cE$ is an isomorphism of abelian groups.\endproof

While Proposition \ref{pres4} describes the left regular representation of the
algebra $\cA_\Z$, Proposition \ref{mapstilderho} allows one to construct a
representation of $\cA_\Z$ on its abelian part $\Z[\Q/\Z]$ as follows.

\begin{prop}\label{repinab} The relations
\begin{equation}\label{abrelations}
\begin{array}{l}
 \theta(x)\,\xi=x\xi \qqq x,\xi \in \Z[\Q/\Z]   \\[3mm]
\theta(\tilde\mu_n)\xi=\tilde\rho_n(\xi)\qqq \xi \in \Z[\Q/\Z],\forall n
 \\[3mm]
\theta(\mu_n^*) \xi =\sigma_n(\xi)\qqq \xi \in \Z[\Q/\Z],\forall n
\end{array}
\end{equation}
define a representation $\theta$ of $\cA_\Z$ on $\Z[\Q/\Z]$.
\end{prop}

\proof It is enough to check that the relations of Definition \ref{overZ} are
fulfilled. The first of  the three relations \eqref{presoverZ1} follows from
\eqref{rhosigmamult}. The second follows from the multiplicativity of
$\sigma_n$. The third one follows again from \eqref{rhosigmamult}. The first
two of the four relations \eqref{presoverZ2} follow from the analogous relation
\eqref{rhomult} on the $\tilde\rho_n$ and $\sigma_n$. The last two relations
both follow from \eqref{divisible}.\endproof

\medskip
\subsection{Relation with the integral Hecke algebra}\label{bchecke}

The original construction of the BC-system \cite{BC} is based on Hecke algebras
of quasi-normal pairs. One considers the inclusion $P_\Z^+ \subset P_\Q^+$
where the ``$ax+b$'' algebraic group $P$  is viewed as the functor which to any
abelian ring $R$ assigns the group $P_R$ of 2 by 2 matrices over $R$ of the
form
\begin{equation}
\label{eq15} P_R = \left\{ \left( \begin{matrix} 1 &b \\ 0 &a
\end{matrix} \right) \, ; \ a,b \in R \, , \ a \ \hbox{invertible}
\right\} \, .
\end{equation}
Here $\Gamma_0=P_\Z^+$ and $\Gamma=P_\Q^+$ denote the restrictions to $a>0$. This inclusion
$\Gamma_0\subset\Gamma$ is such that the orbits of the left action of
$\Gamma_0$ on $\Gamma / \Gamma_0$ are all {\it finite}. The same clearly holds
for orbits of $\Gamma_0$ acting on the right on $\Gamma_0 \backslash \Gamma$.

The integral Hecke algebra ${\mathcal H}_\Z (\Gamma, \Gamma_0)$ is by
definition the convolution algebra of functions of finite support
\begin{equation}\label{fcoset}
f : \Gamma_0 \backslash \Gamma \to \Z ,
\end{equation}
which fulfill the $\Gamma_0$-invariance condition
\begin{equation}\label{finvG0}
f (\gamma \gamma_0) = f(\gamma), \ \ \ \ \forall \gamma \in \Gamma, \forall
\gamma_0 \in \Gamma_0,
\end{equation}
so that $f$ is  defined on $\Gamma_0 \backslash \Gamma / \Gamma_0$. The
convolution product is then defined by the formula
\begin{equation}\label{convhecke}
(f_1 * f_2)(\gamma) = \sum_{\Gamma_0 \backslash \Gamma} f_1 (\gamma
\gamma_1^{-1}) f_2 (\gamma_1) \, .
\end{equation}

There is a presentation of this algebra which is obtained as an extension of
the integral group ring $\Z[\Q/\Z]$  by adjoining elements $\nu_n$ and
$\nu_n^*$ which are formally defined by $\nu_n=\sqrt n \mu_n$, $\nu^*_n=\sqrt n
\mu^*_n$ (with the notations of \cite{BC} \S 4, \ie $\mu_n = n^{-1/2}e_{X_n}$).
The presentation is of the form:
\begin{equation}\label{presoverZ1bis}
\begin{array}{l}
 \nu_n x \nu_n^* = \tilde\rho_n(x)\ \ \ \  \\[3mm]
\nu_n^* x = \sigma_n(x) \nu_n^* \,, \ x  \nu_n =  \nu_n \sigma_n(x),
 \\[3mm]
\nu_{nm}= \nu_n   \nu_m \,, \ \nu_{nm}^* =\nu_{n}^*\nu_{m}^* \qqq n,m\\[3mm]
\nu_n^*  \nu_n  =n \,, \
 \nu_n\nu_m^* = \nu_m^* \nu_n \ \ \ \ (n,m) = 1.
\end{array}
\end{equation}
Comparing this presentation with \eqref{presoverZ1} and \eqref{presoverZ2} one
obtains:

\begin{prop}\label{propmubciso} There exists a unique isomorphism
 \begin{align}\label{mubciso}
 &\phi: ~{\mathcal H}_\Z (\Gamma, \Gamma_0)\stackrel{\sim}{\longrightarrow} \cA_\Z=\Z[\Q/\Z]\rtimes_{\tilde\rho}\N,\notag\\
    &\phi(e(r))=e(r)\qqq r\in \Q/\Z\,, \ \ \phi(\nu_n)=\tilde\mu_n\,, \ \
    \phi(\nu^*_n)=\mu_n^*\,.
\end{align}
\end{prop}

\proof One checks that the relations \eqref{presoverZ1bis} transform into
\eqref{presoverZ1} and \eqref{presoverZ2} under $\phi$.
\endproof

The Hecke algebra ${\mathcal H}_\Z (\Gamma, \Gamma_0)$ admits a natural
involution for which $\nu_n$ and $\nu_n^*$ are adjoint of each other. It is
given (with arbitrary coefficients) by:
\begin{equation}\label{involHecke}
f^* (\gamma) := \overline{f(\gamma^{-1})}, \ \ \ \forall \gamma \in \Gamma_0
\backslash \Gamma / \Gamma_0.
\end{equation}
The rational algebra $\cA_\Z\otimes\Q=\cA_\Q=\Q[\Q/\Z]\rtimes_{\rho}\N$ also
admits a natural involution which coincides with \eqref{involHecke} on the subalgebra
$\Z[\Q/\Z]$ and whose extension to $\cA_\Z$ is dictated by the equation
$\tilde\mu_n=n(\mu_n^*)^*$. \smallskip

Notice that the isomorphism $\phi$ of Proposition \ref{propmubciso} {\it does
not preserve} the involution. The rational subalgebras ${\mathcal H}_\Q
(\Gamma, \Gamma_0)$ and
 $\cA_\Z\otimes_\Z\Q=\cA_\Q$  of the $C^*$-dynamical system
$(C^*(\cG),\sigma_t)$ of Definition \ref{bcsystem} are {\it not the same}. One
has however

\begin{prop}\label{propmubcisoC} The involutive algebras
$\cA_\Z\otimes_\Z\Q=\Q[\Q/\Z]\rtimes_{\tilde\rho}\N$ and ${\mathcal H}_\Q (\Gamma,
\Gamma_0)$ become isomorphic after tensoring by $\C$. An isomorphism is given
by
\begin{equation}\label{mubcisoC}
    \psi(e(r))=e(r)\qqq r\in \Q/\Z\,, \ \ \psi(\nu_n)=n^{-1/2}\,\tilde\mu_n\,, \ \
    \psi(\nu^*_n)=n^{1/2}\,\mu_n^*\,.
\end{equation}
The corresponding  rational subalgebras of the   $C^*$-dynamical system
$(C^*(\cG),\sigma_t)$ are conjugate under $\sigma_{i/2}$.
\end{prop}

\proof As subalgebras of the $C^*$-algebra $C^*(\cG)$, the above involutive
$\Q$-algebras are generated by the $e(r)$ and respectively by  the $\mu_n$ and
$\mu_n^*$ for $\cA_\Z\otimes_\Z\Q$ and by the $\nu_n=n^{1/2}\mu_n$ and
$\nu_n^*=n^{1/2}\,\mu_n^*$ for ${\mathcal H}_\Q (\Gamma, \Gamma_0)$. Thus they
are the same after tensoring with $\C$. To get the conjugacy by $\sigma_{i/2}$,
note that one has $\sigma_t(\mu_n)=n^{it}\mu_n$ and
$\sigma_t(\mu_n^*)=n^{-it}\mu_n^*$. Thus with $\tilde\mu_n=n\mu_n$ one gets
$\sigma_{i/2}(\tilde\mu_n)=n^{1/2}\mu_n=\nu_n$ and
$\sigma_{i/2}(\mu_n^*)=n^{1/2}\mu_n^*=\nu^*_n$.
\endproof

\begin{rem}\label{subtdist} {\rm The above distinction between the two rational
subalgebras of Proposition \ref{propmubcisoC} is overlooked in Proposition 3.25
of \cite{CMbook} Chapter III. However by Proposition \ref{propmubcisoC}, these
two rational algebras are conjugate by $\sigma_{i/2}$ and the $\sigma_t$
invariance of KMS$_\beta$ states thus shows that the values of the restriction
of KMS$_\beta$ states is independent of this distinction. }\end{rem}

\medskip

\section{The endomotive and algebra in characteristic  $p$}\label{sectcharp}

 The group ring $\Z[\Q/\Z]$ together with the endomorphisms $\sigma_n$ and the
 maps $\tilde\rho_n$ give a model over $\Z$ of the BC-endomotive.

 In this section we study the reduction of this model at a prime $p$
 both at the level of the endomotive and  of the noncommutative
 crossed product algebra.
 From now and throughout this section we shall work over a perfect field $\K$ of
characteristic $p>0$, such as a finite extension $\F_q$ of $\F_p$ or a
separable closure of $\F_p$.

 We first show that, by applying reduction at $p$ and specializing $n$ to be
 $p^\ell$, the endomorphism $\sigma_n$ on $\K[\Q/\Z]$ is identified with the geometric
 Frobenius correspondence. The group
 algebra   $\K[\Q/\Z]$ decomposes as a tensor product of the group algebra
 $\K[\Q_p/\Z_p]$  of
 the
 $p$-torsion $\Q_p/\Z_p$ of $\Q/\Z$ by the group algebra of fractions with denominators
 prime to $p$. The structure of the latter algebra is essentially insensitive
 to characteristic $p$. The new fact specific to characteristic $p$ is that
 the group algebra $\K[\Q_p/\Z_p]$ is unreduced and in fact local. We
 concentrate on this ``$p$-part" of the abelian algebra.

\smallskip

We then form a new noncommutative  algebra obtained  as the crossed  product of
the $p$-part $\K[\Q_p/\Z_p]$ by the sub-semigroup of $\N$ given by powers of
$p$. We exhibit the nilpotent nature of this algebra by showing that it admits
a faithful representation as infinite triangular matrices.

\subsection{The endomotive in characteristic $p$}

  The relevant properties of the algebra $\cA_\Z(\K)$ can be isolated
by decomposing the groups as follows
\begin{equation}\label{groupdec}
\Q/\Z=\Q_p/\Z_p\times (\Q/\Z)^{(p)} \,.
\end{equation}
Here $\Q_p/\Z_p$ is identified with the subgroup of $\Q/\Z$ of fractions with
denominator a power of $p$ and  $(\Q/\Z)^{(p)}$ is interpreted as the subgroup
of $\Q/\Z$ of fractions with denominator prime to $p$. At the group algebra
level one gets
\begin{equation}\label{galevel}
   \K[\Q/\Z]=\K[\Q_p/\Z_p]\otimes \K[(\Q/\Z)^{(p)}]\,.
\end{equation}

The  decomposition
\begin{equation}\label{semigrdec}
\Q_+^*=p^{\Z}\times \Q_+^{(p)}
\end{equation}
 corresponds to the decomposition of the
semigroup $\N$ as a product of the semigroup $p^{\Z^+}$ of powers of $p$ and
the semigroup $\N^{(p)}$ of numbers prime to $p$.
  There is no essential difference with the characteristic zero set-up for the action of
   $\N^{(p)}$  on $\K[\Q/\Z]$.  In fact,
    the endomorphism $\rho_n$ on $\K[\Q/\Z]$ retains a meaning when $n$
    is prime to $p$, since the denominators  in the definition of the projection
$\pi_n$ and of $\rho_n$ (\ie the partial inverse of $\sigma_n$) are integers
prime to p. Moreover, we notice that when $n$ is prime to $p$ the equation
$T^n-1=0$ is unramified. On the other hand, when $n\in p^{\Z^+}$ there is no
way to single out the component of $\{1\}$ in $\Sp(A_n)$ since in that case the
above equation has $1$ as a multiple root. Therefore, our study will focus on
the understanding of\medskip

$\bullet$~The endomorphism $\sigma_n$, for $n=p^\ell\in p^{\Z^+}$\smallskip

$\bullet$~The algebra  $A_{p^\infty}=\varinjlim_\ell
A_{p^\ell}\otimes_\Z \K=\K[\Q_p/\Z_p]$\medskip

We first  show the relation between the endomorphisms $\sigma_n$, for
$n=p^\ell\in p^{\Z^+}$, and the (relative) geometric Frobenius homomorphism
acting on the algebra $A_{p^\infty}$.

\begin{prop}\label{geomfrobprop}
Let $\sigma_{\F_p}\in \Aut( \K)$ be the small Frobenius automorphism given by
$\sigma_{\F_p}(x)=x^p$ for all $x\in \K$, then for any $\ell$,
\begin{equation}\label{sigptaup1}
    (\sigma_{p^\ell}\otimes \sigma_{\F_p}^\ell)(f)=f^{p^\ell}\qqq f\in
     \K[\Q/\Z]=\Z[\Q/\Z]\otimes_\Z \K,\,\ \ \ \forall \ell\in\N.
\end{equation}
\end{prop}

\proof Both sides of \eqref{sigptaup1} define an endomorphism of the ring, thus
it is enough to check that they agree on elements of the form $e(r)\otimes x$,
for $x\in\K$ and $r\in \Q/\Z$. One has
$$
(\sigma_{p^\ell}\otimes \sigma_{\F_p}^\ell)(e(r)\otimes x)=e(p^\ell r)\otimes
x^{p^\ell}=(e(r)\otimes x)^{p^\ell}\,,
$$
which gives the required equality.
\endproof

\begin{rem}\label{geomFr}
The relation $(\sigma_{p^\ell}\otimes \sigma_{\F_p}^\ell)(f)=f^{p^\ell}$ of
\eqref{sigptaup1} shows that we can interpret $\sigma_{p^\ell}$ as the
Frobenius correspondence acting on the pro-variety $(\mu^{(\infty)}\otimes_\Z
\K)$. This means that at the fixed level $\mu^{(m)} = \text{Spec}(A_{m})$, the
definition of $\sigma_{p^\ell}$ coincides with the Frobenius morphism $\varphi$
as in \cite{Tate}, p. 24 (\ie the morphism inducing in \'etale cohomology the
geometric Frobenius $\Phi$ of Deligne).
\end{rem}

An equivalent description of the algebra $A_{p^\infty}=\varinjlim_\ell
A_{p^\ell}\otimes_\Z \K=\K[\Q_p/\Z_p]$ will be given in terms of the following
(local) convolution algebra of functions which displays an explicit choice of a
basis. At a fixed   level $\ell$, \ie for the algebra $A_{p^\ell}$, this choice
of basis corresponds to the powers $\epsilon^k$ of the natural generator
$\epsilon= \delta_{p^{-\ell}}$, as in \eqref{deltaadefn},  of the maximal ideal
of the local ring $A_{p^\ell}$, \cf Remark \ref{filtration}.

\begin{defn} We define the algebra $\cT(p)$ (over $\K$) as the convolution algebra of
$\K$-valued functions with finite support on the semigroup $S_+=\cup
\frac{1}{p^n}\Z^+\subset \R$, modulo the ideal of functions with support in
$S\cap [1,\infty[$, with the convolution product given by
\begin{equation}\label{conv}
    f_1\star f_2(c)=\sum_{a+b=c}f_1(a)f_2(b).
\end{equation}
\end{defn}

We extend any function $f$ on $S_+$ to a function on $S=\cup \frac{1}{p^n}\Z$
which fulfills $f(a)=0$ for all $a<0$. This extension is compatible with the
convolution product. By construction the algebra $\cT(p)$ is commutative and
{\em local}. It has a unique character: the homomorphism of $\K$-algebras given
by evaluation at $0$, that is,
\[
\epsilon_0: \cT(p)\to \K\,, \ \ f\mapsto f(0).
\]
Any element in the kernel of this character is nilpotent. The kernel
$\text{Ker}(\epsilon_0)$ of this character is the only maximal ideal.

For any $a\in S\cap [0,1)$ we let $\delta_a\in \cT(p)$ be given by
\begin{equation}\label{deltaadefn}
  \delta_a(a)=1\,, \ \   \delta_a(b)=0 \ \text{if} \ b\neq a\,.
\end{equation}
Any $f\in \cT(p)$ is a finite sum $f=\sum f(a) \delta_a$ and $\delta_0$ is the
unit ${\bf 1}$ of the algebra $\cT(p)$.

\smallskip

\begin{prop}\label{triangular} 1) The following map induces in the limit an isomorphism of
$A_{p^\infty}$ with $\cT(p)$,
\begin{equation}\label{triangmap}
  \iota\;:\;  A_{p^\ell}\otimes_\Z \K\to \cT(p)\,, \ \
  \iota(e(1/p^\ell))= {\bf 1}-\delta_{p^{-\ell}}.
\end{equation}
2) The endomorphism $\sigma_p$ corresponds by the above isomorphism to the
following endomorphism of $\cT(p)$,
\begin{equation}\label{sigptaup}
    \sigma_p(f)(a)=f(a/p)\qqq f\in \cT(p).
\end{equation}

3) The map $\tilde\rho_p$ corresponds by the above isomorphism to the following
map of $\cT(p)$,
\begin{equation}\label{rhoptaup}
    \tilde\rho_p(f)(a)=f(p\,a-(p-1))\qqq f\in \cT(p).
\end{equation}
\end{prop}

\begin{figure}
\begin{center}
\includegraphics[scale=0.5]{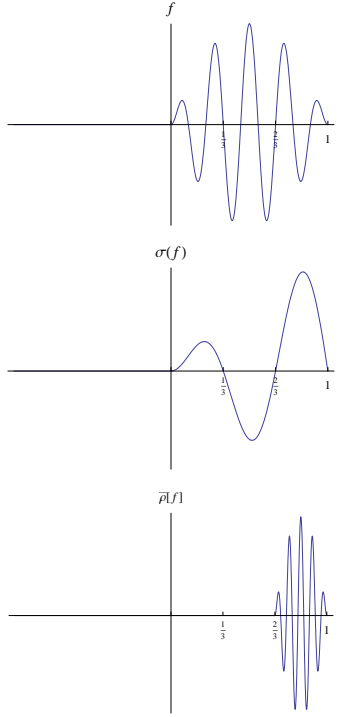}
\end{center}
\caption{The maps $\sigma_p$ and $\tilde\rho_p$ \label{threegraphsfig} }
\end{figure}

Note that both maps $\sigma_p$ and $\rho_p$ are given by an affine change of
variables as shown in Figure \ref{threegraphsfig}.

\proof 1) Let us check that the $\iota(e(1/p^\ell))$ fulfill the rules of the
generators $e(1/p^\ell)$. In characteristic $p$ one has
$$
({\bf 1}-T)^{p^\ell}={\bf 1}-T^{p^\ell}
$$
Thus to show that $\iota(e(1/p^\ell))^{p^\ell}={\bf 1}$ it is enough to check
that the $p^\ell$ power of the characteristic  function $\delta_{p^{-\ell}}$ is
equal to $0$. This follows from the equalities $\delta_a\star
\delta_b=\delta_{a+b}$ (using \eqref{conv}) and $\delta_1=0$. In fact one needs
to show that
$$\iota(e(1/p^\ell))^p=\iota(e(1/p^{\ell-1}))$$
which means that
$$
({\bf 1}-\delta_{p^{-\ell}})^p={\bf 1}-\delta_{p^{-\ell+1}}
$$
and this follows from $\delta_{p^{-\ell}}^p=\delta_{p^{-\ell+1}}$.

2) It is enough to check \eqref{sigptaup} on the elements
$\iota(e(1/p^\ell))={\bf 1}-\delta_{p^{-\ell}}=\delta_0-\delta_{p^{-\ell}}$.
The right hand-side of \eqref{sigptaup} defines an endomorphism of $\cT(p)$
which transforms $\delta_a$ into $\delta_{pa}$ and this gives
$\sigma_p(e(1/p^\ell))=e(1/p^{\ell-1})$ as required.

3) Note that since $f$ is extended to a function on $S=\cup \frac{1}{p^n}\Z$
which fulfills $f(a)=0$ for all $a<0$, the formula \eqref{rhoptaup} makes sense
and the function $\rho_p(f)$ vanishes on the interval $[0,\frac{p-1}{p})$. In
characteristic $p$ one has, with $q=p^\ell$,
\begin{equation}\label{charpequ}
    \sum_0^{q-1}\, T^k=({\bf 1}-T)^{q-1}
\end{equation}
since multiplying both sides by $({\bf 1}-T)$ gives $({\bf 1}-T)^{q}$. This
shows that
\begin{equation}\label{charpequ1}
\iota(\tilde\rho_p(1))=\sum_0^{p-1}\iota(e(k/p))=({\bf
1}-\iota(e(1/p))^{p-1}=\delta_{\frac{p-1}{p}}\,.
\end{equation}
Now, by \eqref{rhosigmamult}, one has
$$
\tilde\rho_p(\sigma_p(x))=x  \tilde\rho_p(1)\qqq x\in A_{p^\infty}
$$
which gives the required equality \eqref{rhoptaup} using the surjectivity of
$\sigma_p$ and the fact that in the algebra $\cT(p)$ the convolution by
$\iota(\tilde\rho_p(1))=\delta_{\frac{p-1}{p}}$ is given by the translation by
$\frac{p-1}{p}$.
\endproof

\begin{cor}\label{corppart} The kernel of $\sigma_p$ is the nilpotent ideal
\begin{equation}\label{nilideal}
    \ker\, \sigma_p=\{f \in \cT(p)\,|\, f(a)=0\qqq a\in [0,\frac 1p)\}
\end{equation}
\end{cor}

\proof One has $f\in \ker\, \sigma_p$ iff $f(a/p)=0$ for all $a\in S\cap [0,1)$
which gives \eqref{nilideal}. Any element $f\in \ker\, \sigma_p$ thus fulfills
$f^p=0$ in the algebra $\cT(p)$.
\endproof\bigskip

\subsection{The BC-algebra in characteristic $p$}

By definition, the BC-algebra over $\K$ is the tensor product:
\begin{equation}\label{bcoverk}
    \cA_\Z(\K)=\cA_\Z\otimes_\Z\K
\end{equation}
 By Corollary \ref{corpres4} the
$\K$-linear space $\cA_\Z(\K)$ coincides with the vector space $\K[\Q/\Z\times
\Q_+^*]$ and because of that we will work with the corresponding linear basis
of monomials \eqref{basismonZbis}.

The remaining part of this section is dedicated to the study of the algebra
\begin{equation}\label{ppowers}
   \cC_p= A_{p^\infty}\rtimes_{\tilde\rho}p^\N \sim \cT(p)\rtimes_{\tilde\rho}p^\N
\end{equation}
We shall refer to $\cC_p$ as to the $p$-part of the algebra $\cA_\Z(\K)$. We
keep the same notation as in \S~\ref{bczz}.

\begin{lem}\label{pres1} The following monomials form a linear basis of the
algebra $\cC_p$:
\begin{equation}\label{basismon}
    \tilde\mu_p^n\,\delta_a \,, \ n\in \N\,, \ a\in S\cap [0,1) \,, \ \
    \delta_a {\mu_p^*}^m \,, \ m\geq 0\,, \ a\in S\cap [0,1)\,.
\end{equation}
\end{lem}

\proof For $m=0$ we use the notation ${\mu_p^*}^0=1$ so that the above
monomials contain the algebra $A_{p^\infty}$ and the generators $\tilde\mu_p$
and $\mu_p^*$. Thus it is enough to show that the linear span of these
monomials is stable under the product. One has
\begin{equation}\label{basismon1}
\tilde\mu_p^n\,x\,\tilde\mu_p^m\,y=\tilde\mu_p^{n+m}\sigma_p^m(x)y\,,\ \ \ x
{\mu_p^*}^n\,y {\mu_p^*}^m=x\sigma_p^n(y){\mu_p^*}^{n+m}
\end{equation}
and, for $n>0$,
\begin{equation}\label{basismon2}
x{\mu_p^*}^n\,\tilde\mu_p^m\,y=0
\end{equation}
while
\begin{equation}\label{basismon3}
\tilde\mu_p^n\,x\,y {\mu_p^*}^m=\begin{cases} \tilde\rho_p^n(xy){\mu_p^*}^{m-n}&\text{if}~m\geq n\\
\tilde\mu_p^{n-m}\tilde\rho_p^m(xy) &\text{if}~n>m.\end{cases}
\end{equation}
which shows that the linear span of the above monomials is an algebra.
\endproof

Note that Lemma \ref{pres1} also follows directly from Lemma \ref{pres1Z}.\smallskip

In order to exhibit the {\em nilpotent} nature of this algebra we  now show
that the representation of Proposition \ref{repinab} is given by infinite
triangular matrices.

 We let $\K[S\cap[0,1)]$ be the $\K$-linear space of
$\K$-valued functions with finite support on $S\cap[0,1)$ and denote by
$\xi_a$, $a\in
 S\cap[0,1)$ its canonical basis. For $a\in S_+, a\geq 1$ we let $\xi_a=0$ by
 convention.

 Let $G=S\rtimes \Z$ be the
semi-direct product of the additive group $S$ by the action of $\Z$ whose
generator acts on $S$ by multiplication by $p$. The group $G$ acts on $S$ by
affine transformations
\begin{equation}\label{gsemigroup}
\alpha_{g}(b)=p^n b+ a\qqq g=(n,a)\in G\,, \ b\in S\,.
\end{equation}

\begin{lem}\label{defgsemigroup} Let $G$ be the group defined above.
\begin{enumerate}
  \item The condition
  \begin{equation}\label{trisemi}
    x\geq a\Rightarrow g(x)\geq a\qqq a\in [0,1]
  \end{equation}
defines a sub-semigroup $G^+\subset G$.
  \item $G^+$ acts on $\K[S\cap[0,1)]$ by
  \[
  \tau(g)\xi_a=\xi_{g(a)}.
  \]
  \item The semi-group $G^+$ is generated by the elements $g_a=(0,a)$, for $a\in
  S\cap[0,1)$, $\alpha=(-1, (p-1)/p)$, $\beta=(1,0)$.
\end{enumerate}
\end{lem}

\proof The first statement is obvious. The second follows from \eqref{trisemi}
for $a=1$. Let us prove (3). Let $g=(n,a)\in G^+$. If $n=0$ then $g=g_a$ with
$a\geq 0$. If $n>0$ then $g(x)=p^nx+a$ and taking $x=0$ shows that $a\geq 0$ so
that $g=g_a\beta^n$. For $n=-m<0$, $g(x)=p^{-m}x+a$. Taking $x=1$ and using
\eqref{trisemi}, one gets $p^{-m}+a\geq 1$ \ie $a=1-p^{-m}+b$ for some $b\geq
0$. Thus $g=g_b\alpha^m$.\endproof

\smallskip
\begin{prop} \label{proptriangular}  The   equations
\begin{equation}\label{triangact}
\begin{array}{l}
 \theta(\delta_b)\,\xi_a=\tau(g_b)\xi_a=\xi_{a+b}\ \ \ \  \\[3mm]
\theta(\tilde\mu_p)\,  \xi_a= \tau(\alpha)\xi_a= \xi_{\frac{a+p-1}{p}},
 \\[3mm]
\theta(\mu_p^*)\,  \xi_a =\tau(\beta)\xi_a=\xi_{pa}
\end{array}
\end{equation}
define a faithful representation
\[
\cC_p \stackrel{\theta}{\to} \text{End}(\K[S\cap[0,1)])
\]
 of
 the algebra $\cC_p$ by (lower)-triangular matrices $T=(T_{a,b})$, $T_{a,b}\in \K$ with $a,b\in
 S\cap[0,1)$.

\end{prop}

\proof The matrix associated to  $T\in \text{End}(\K[S\cap[0,1)])$  is
defined by
\begin{equation}\label{matrix}
(Tf)(a)=\sum T_{a,b}\, f(b)\qqq f\in \K[S\cap[0,1)]
\end{equation}
Thus the matrices associated to the operators given in \eqref{triangact} are:
\begin{equation}\label{triangactmatrix}
\begin{array}{l}
 (\delta_a)_{c,d}=0\,\ \text{if}\ c-d\neq a \,,\ \ (\delta_a)_{c,c-a}=1  \\[3mm]
(\tilde\mu_p)_{c,d}=0\,\ \text{if}\ c\neq \frac{d+p-1}{p}\,,\  \
(\tilde\mu_p)_{c,d}=1\,\ \text{if}\ c= \frac{d+p-1}{p}
 \\[3mm]
(\mu_p^*)_{c,d}=0\,\ \text{if}\ c\neq pd\,, \ (\mu_p^*)_{c,c/p}=1\,.
\end{array}
\end{equation}
They are lower triangular. Indeed one has
$$
c\geq c-a\,, \ \ \frac{d+p-1}{p}>d\,, \ \ p\,d\geq d  \qqq c,d\in
 S\cap[0,1)\,.
$$
  One then needs to show that the
defining relations of the algebra $\cC_p$ are fulfilled. These relations are
obtained from the presentation of Definition \ref{overZ} by restriction to the
$p$-part. Thus they are fulfilled by specializing Proposition \ref{repinab} to
the $p$-part. One can also check them directly. By construction the action of
the $\delta_a$ gives a representation of the convolution algebra $\cT(p)$.
 The three
relations of \eqref{presoverZ1}
\begin{equation}\label{presoverZ1bis2}
\begin{array}{l}
\tilde\mu_p\, x \mu_p^* = \tilde\rho_p(x)\ \ \ \  \\[3mm]
\mu_p^*\, x = \sigma_p(x) \mu_p^*  \\[3mm]
x \tilde\mu_p = \tilde\mu_p \sigma_p(x),
\end{array}
\end{equation}
follow directly from the group action $\tau$.  Moreover one has the additional
relation
$$
\mu_p^*\,\tilde\mu_p\, \xi_a=\mu_p^*\,\xi_{\frac{a+p-1}{p}}=\xi_{ a+p-1 }=0,
$$
which corresponds to the third relation of \eqref{presoverZ2}. Its validity
follows from  $\xi_b=0$ for all $b\geq 1$.

Let us prove that the representation is faithful. By Lemma \ref{pres1} any
element  $x\in\cC_p$ is a finite linear combination $x=\sum \lambda_j
\tau(g_j)$ ($\lambda_j\in \K$) of monomials with  $g_j\in G^+$.  Now for any
two distinct elements $g, h\in G^+$ the set of elements $a\in S$ such that
$g(a)=h(a)$ contains at most one element. Thus, since $g_j(0)\in [0,1)$, one
can find an element $b\in S\cap [0,1)$ such that
$$
g_j(b)\in S\cap [0,1)\qqq j \,, \ \ g_j(b)\neq g_k(b)\qqq j\neq k\,.
$$
We then have
$$
x\,\xi_b=\sum \lambda_j \tau(g_j)\,\xi_b=\sum \lambda_j \xi_{g_j(b)}
$$
and $x\xi_b\neq 0$ if $x\neq 0$. Thus the representation is faithful.
\endproof

\begin{rem}\label{filtration} {\rm By construction the algebra $A_{p^\infty}$ is the inductive limit of
 the local rings $A_q=\K[T]/(T^q-1)$, $q=p^\ell$. We let
 \begin{equation}\label{tildepi}
    \tilde\pi_q =\tilde\rho_q(1)= 1+T+\cdots +T^{q-1}\in A_q=\Z[T]/(T^q-1).
\end{equation}
The local ring $\K[T]/(T^q-1)$ is generated over $\K$ by the nilpotent element
$\epsilon = T-1$ ($\epsilon^q = 0$). The principal ideal of multiples of
$\epsilon$ is the maximal ideal. We use the natural decreasing finite
filtration by powers of the maximal ideal,
\begin{equation}\label{filtr}
    F^i(\K[T]/(T^q-1))=\epsilon^i\,\K[T]/(T^q-1)
\end{equation}
One has, using \eqref{charpequ1},
\begin{equation}\label{filtr1}
\tilde\pi_q =\epsilon^{q-1}\in F^{q-1}(\K[T]/(T^q-1))\,, \ \
F^q(\K[T]/(T^q-1))=\{0\}
\end{equation}
 Thus the operator $\tilde\pi_q$ detects the top piece
of the filtration.

\smallskip
The following equalities show that the subalgebra $\cP\subset A_{p^\infty}$
generated by the  $\tilde\rho_p^m(1)=\tilde \pi_{p^m}=\tau_m$ is stable under
the $\tilde\rho_p$ and $\sigma_p$.
\begin{equation}\label{basicpropeps}
   \tau_m\tau_n=0\,, \ \ \tilde\rho_p^m(\tau_n)=\tau_{m+n}
   \,, \ \ \sigma_p(\tau_n)=0 \qqq m,n\in \N\,.
\end{equation}
As above one checks that  the following monomials form a linear basis of the
crossed product algebra $\cP\rtimes_{\tilde\rho}p^\N$:
\begin{equation}\label{basismonp}
    \tilde\mu_p^n\,\tau_k \,, \ n\in \N\,, \ k\geq 0 \,, \ \
    \tau_k  {\mu_p^*}^m \,, \ m\geq 0\,, \ k\geq 0\,.
\end{equation}
Since $\cP\rtimes_{\tilde\rho}p^\N$ is a subalgebra of the algebra $\cC_p$,
Proposition \ref{proptriangular} yields in particular, a triangular representation of
$\cP\rtimes_{\tilde\rho}p^\N$.

 }\end{rem}

\subsection{The effect of reduction}\label{sectreduc}

In the original definition of endomotives given in \cite{CCM}, we assumed that
the algebras are reduced. This is in agreement with the classical definition of
Artin motives (\cf~\cite{LNM900}, II p.~211). In the present context, namely
working over a perfect field $\K$ of characteristic $p$, one can still restrict
to reduced algebras by functorial reduction. One can see in the result below
that this reduction introduces a drastic simplification of the algebra, which,
in particular, eliminates the problem of denominators.

\begin{prop}\label{triangularbis} The reduced algebra of $\varinjlim_n
A_n\otimes_\Z \K$ is the group ring over $\K$ of the subgroup of $\Q/\Z$ of
fractions with denominator prime to $p$. Moreover, $\sigma_p$ induces an
automorphism on the reduced algebra.
\end{prop}

\proof This amounts to showing that for $n=p^km$ with $m$ prime to $p$, the
reduced algebra of $\K[T]/(T^n-1)$ is the algebra $\K[T]/(T^m-1)$. If $n=p^km$
then the group $\Z/n\Z$ splits canonically as a product of $\Z/p^k\Z$ and
$\Z/m\Z$. At the group ring level, this corresponds to a tensor product
decomposition. Since the reduction in characteristic $p$ of the group ring of
$\Z/p^k\Z$ is the ground field $\K$, the first factor in the tensor product
disappears and the reduction only leaves the second factor. This proves the
first statement. It is then enough to observe that, for $m$ prime to $p$, the
multiplication by $p$ is an automorphism of $\Z/m\Z$. Since $\sigma_p$
preserves the levels, this is compatible with the map of the inductive system
of algebras.
\endproof

\begin{cor}\label{nodenom}
In the case of the reduced algebra in characteristic $p$, the inverses $\rho_n$
of the $\sigma_n$ only involve denominators that are prime to $p$.
\end{cor}

\proof The case where $n$ is prime to $p$ is clear. Suppose that $n=p^k$. Then
by Proposition \ref{triangularbis}, $\sigma_n$ is an automorphism of the
reduced algebra  since multiplication by $n$ is an automorphism of the group
$(\Q/\Z)^{(p)}$. One then defines $\rho_n$ as its inverse and the corresponding
$\pi_n$ is then equal to one, since $\sigma_n$ is injective.
\endproof

Note that passing from $\Q/\Z$ to the subgroup $(\Q/\Z)^{(p)}$ (\ie the
prime-to-$p$ component) is the same, when dealing with the Pontrjagin dual
groups, as removing from the ring of finite adeles the component at $p$. This
suggests that there is a connection with the localized system at $p$ in
\cite{CCM2} (\cf~Definition~8.14 and Theorem~8.15). Note however that unlike
the setting of \cite{CCM2}, here the coefficients are taken in a field of
positive characteristic, so that the notion of KMS states should be taken in
the extended sense of \cite{CM}.

\begin{rem}\label{confuse}{\rm Notice that reducing the abelian part of the
algebra and then taking the crossed product as we did in this section is not
the same thing as moding out the crossed product algebra $\cA_\Z$ by its
nilpotent radical. }\end{rem}

 \subsection{Endomotives in the unreduced case}\label{sectendounred}

As we have seen in the previous sections, when taking coefficients in a field
of positive characteristic the
  BC-endomotive
  involves unreduced finite dimensional commutative algebras which strictly
  speaking do not correspond to classical Artin motives.
The construction in characteristic $p$ that we gave in the case of the BC-algebra in fact extends to a more general class of endomotives constructed
 from finite, self maps of algebraic varieties as in \cite{CCM}, but without requiring that these maps are unramified over the base point.

This leads us naturally to consider the problem of a general construction of
endomotives in arbitrary characteristic.   Roughly speaking, an
  endomotive over a (perfect) field $\K$ is given by assigning:\smallskip

  $\bullet$~An inductive system of augmented
    commutative $\K$-algebras, finite dimensional as vector spaces over
  a perfect field $\K$ (\ie Artinian commutative $\K$-algebras).\smallskip

  $\bullet$~A commutative family of correspondences $\sigma_n$.\smallskip

This set of data should of course be compatible with the constructions that we have developed in this paper as well as in \cite{CCM}, namely
\smallskip

1)~It should determine homomorphisms (correspondences) such as the $\rho_n$'s, when denominators (\ie division by $n$) are allowed.\smallskip

2)~It should be fulfilled by the endomotives associated to self-maps of pointed varieties as described in \cite{CCM} (Example~3.4).\smallskip

\section{The BC endomotive over $\F_1$}\label{BCF1sec}

In this section we show that the BC endomotive has a model defined over $\F_1$
from which one recovers the original endomotive by extension of scalars to $\Q$.

 \begin{prop}\label{propdeffun}
\begin{enumerate}
  \item[a)] The BC-endomotive has a model over $\F_1$.
    \item[b)] The original BC-endomotive is obtained by extension of the
  scalars from $\F_1$ to $\Q$.
  \end{enumerate}
  \end{prop}

  \proof We start with the projective system of affine varieties $\mu^{(n)}$ over $\F_1$, defined as in \S~\ref{F1roots1} and Proposition \ref{mukF1}. This system shows that the abelian part of the BC endomotive is defined over $\F_1$. Notice that these are pointed varieties because the algebras $A_n$ are naturally augmented. The augmentations fit together in the inductive system of algebras because they come from the natural augmentation of the group ring.

It remains to show that the $\sigma_n$ are morphisms of varieties over $\F_1$
in the sense of \cite{Soule}. This is a consequence of the construction of the
projective system of the varieties $\mu^{(n)}$ over $\F_1$, as these are
obtained by applying the functor $\cF$ from varieties over $\Z$ to gadgets over
$\F_1$ (\cf~Proposition \ref{mukF1}). Notice that the maps $\sigma_n$ preserve
levels and are given at each level $A_n$ by \eqref{actionlevelk}. Thus, the
$\sigma_n$ are morphisms in the category of   varieties over $\Z$, and as such
they define morphisms of varieties over $\F_1$ through the functor $\cF$.
\endproof

\subsection{The automorphisms of $\F_{1^\infty}/\F_1$ and the symmetries of the BC system}\label{frobf1}

In \cite{Haran}, the analog of the Frobenius automorphism for the extension
$\F_{1^\infty}$ of $\F_1$ is described as follows. Suppose given a set $X$ with
a free action of the roots of unity (that is a vector space over $\F_{1^\infty}$
when one adds an extra fixed point $0$). Then, given an element $\alpha \in
\hat\Z^*$ (and more in general a non-invertible one in $\hat\Z$) one defines a
new action on the same set by the rule
\begin{equation}\label{xaction}
\zeta: x \mapsto \zeta^\alpha x.
\end{equation}
For $\alpha = n$ an integer, this means that one replaces the action of a root of unity
$\zeta$ by that of $\zeta^n$. The $\hat\Z$-powers of the Frobenius are then
defined by setting
\begin{equation}\label{psitoalpha}
\psi^\alpha:\F_{1^\infty} \to \F_{1^\infty}
\end{equation}
to be the map that sends the action of roots of unity $\zeta$ on a given
$\F_{1^\infty}$-vector space $X$ to the action by $\zeta^\alpha$.

When reformulated additively, after making an identification of the group of roots of unity
with $\Q/\Z$, one can write the action \eqref{xaction} in the form
\begin{equation}\label{xaction2}
e(r) \mapsto e(\alpha(r)),
\end{equation}
where $\alpha\in\hat\Z$ is seen as an element in $\hat\Z=\Hom(\Q/\Z,\Q/\Z)$.

Thus, in terms of the BC system, the Frobenius appears naturally in the semigroup action.
This is the case of \eqref{xaction2} where $n$ is an integer and it gives the
action of  the $\sigma_n$. Moreover, it is also important to keep in mind that the Frobenius action \eqref{xaction2} also
recovers the symmetries of the BC system. In fact, the symmetries by
automorphisms, given by $\hat\Z^*$ act exactly like the corresponding Frobenius
$\psi^{\hat\Z^*}$ (\cf \cite{BC}, \cite{CCM}).

\subsection{The Frobenius correspondence and the BC endomotive over $\F_1$}

We now show that not only the BC-endomotive has a model over $\F_1$, but in
fact it captures the structure of the   extension $\F_{1^\infty}=\varinjlim
\F_{1^n}$ over $\F_1$ by means of the Frobenius correspondence.

\begin{thm}\label{FrobBCF1}
The structure of the BC-endomotive   corresponds to the structure of $\F_{1^\infty}$ over $\F_1$ as follows\smallskip

a)~The abelian part of the BC-endomotive over $\F_1$ corresponds
to the inductive system of ``extensions'' $\F_{1^n}$.\smallskip

b)~The endomorphisms $\sigma_n$ describe the Frobenius correspondence, in the
sense that on the algebra $\Z[\Q/\Z]\otimes_\Z \K$, for $\K$ a perfect field of
characteristic $p>0$, the endomorphisms $\sigma_n$, $n=p^\ell$ ($\ell\in\N$)
coincide with the  Frobenius correspondence described in Remark~\ref{geomFr}.
\end{thm}

\proof For a), we recall that the abelian part of the BC-endomotive
over $\F_1$ is defined by the projective system $\mu^{(\infty)}$
of algebraic varieties $\mu^{(n)} = \cF\text{Spec}(A_n)$, $A_n = \Z[\Z/(n)]$
\cf\eqref{defnmun}. By means of the isomorphism of algebras
$\Z[T]/(T^n-1) \stackrel{\sim}{\to} A_n$, $u(n)\mapsto e(\frac{1}{n})$
the inductive system of extensions $\F_{1^n}\subset\F_{1^m}$ ($n|m$)
corresponds, after extending the coefficients to $\Z$,
to the projective limit $\cdots\to\mu^{(m)}\to \mu^{(n)}\to$
which defines geometrically the abelian part of the BC-endomotive.
We refer to \S~\ref{F1roots1} for the details.

For b), we refer to  Proposition~\ref{geomfrobprop}. We recall that on an
algebraic variety $X_0$ defined over a finite field $\F_q$ ($q=p^\ell$) the
Frobenius morphism $\varphi: X_0 \to X_0$ satisfies the property that the
composition $\varphi\times\sigma: \bar X_0 \to \bar X_0$, where $\sigma:
\bar\F_q \to \bar\F_q$ is the arithmetic Frobenius automorphism, acts on $\bar
X_0 := X_0\times \bar\F_q$ by fixing points and by mapping $f\mapsto f^q$ in
the structure sheaf of $X_0$. Here,
 $f^q$ denotes the section $f$ whose coefficients are raised
 to the $q$-th power. At each fixed level $A_{n}$ of
 the inductive system of algebras $A_n\otimes_\Z \K$, the endomorphisms
 $\sigma_k$, for $k = p^\ell$, behave in exactly
 the same way as the Frobenius homomorphisms
 (\cf~Proposition~\ref{geomfrobprop}).
\endproof

\subsection{Recovering the analytic endomotive}\label{analyendo}

In \cite{Soule}, the set of data which define a variety $X$ (of finite type)
over $\F_1$ is inclusive of the important {\em analytic} information supplied
by the assignment of a commutative Banach $\C$-algebra $\cA_X$
(\cf\S~\ref{F1roots1} of this paper). The definition of $X$ implies that
functions of $\cA_X$ can be evaluated at the points of $X$.  We shall now show
that this analytic part of the set of data which define the BC-system as a
pro-variety over $\F_1$ supply naturally the structure of an analytic
endomotive in the sense of \cite{CCM}. The point is that the set-up which
describes the pro-variety $\mu^{(\infty)}$ is inclusive of
 the information supplied by an {\em inductive} system
 of Banach $\C$-algebras $\cA_{\mu^{(n)}}=A_n\otimes_\Z\C$,
 \cf\eqref{trucV}. Taking the inductive limit of these yields the algebra
\begin{equation}\label{indlimA}
    \cA_{\mu^{(\infty)}}=\varinjlim \cA_{\mu^{(n)}}
\end{equation}
since the functor  $X\mapsto\cA_X$ is contravariant. The following statement is
a direct consequence of the construction of the model of the BC-endomotive over
$\F_1$ and of \eqref{indlimA}:

\begin{prop}\label{endoBCanal}
The analytic part of the pro-variety over $\F_1$ associated to the BC
endomotive over $\F_1$ coincides with the analytic endomotive of the BC system
as described  in \cite{CCM}.
\end{prop}

\end{document}